\input jytex.tex   % available from hep-th
\typesize=10pt \magnification=1200 \baselineskip17truept
%\baselineskip25truept
\footnotenumstyle{arabic} \hsize=6truein\vsize=8.5truein
%\input Spinharm.lab
%\draft
%\leftmargin=1.25in
%\oddleftmargin=.5in
%\evenleftmargin=1.5in
\sectionnumstyle{blank}
\chapternumstyle{blank}
\chapternum=1
\sectionnum=1
\pagenum=0
%\referencestyle{preordered}
% title style follows

\def\begintitle{\pagenumstyle{blank}\parindent=0pt
\begin{narrow}[0.4in]}
\def\endtitle{\end{narrow}\newpage\pagenumstyle{arabic}}

% exercise style follows

\def\beginexercise{\vskip 20truept\parindent=0pt\begin{narrow}[10
truept]}
\def\endexercise{\vskip 10truept\end{narrow}}

% **************    my jyTeX abbreviations   *****************

\def\eql#1{\eqno\eqnlabel{#1}}
\def\ref{\reference}
\def\peq{\puteqn}
\def\pref{\putref}

\def\mgn{\marginnote}
\def\bex{\begin{exercise}}
\def\eex{\end{exercise}}

% *********************** My definitions ************************

\font\open=msbm10 %scaled\magstep1 % For VAX. Borde p195.

 %scaled\magstep1 % For VAX. Borde p195.
%\font\open=msym10 %scaled\magstep1 % For Arbortxt on PC
%\font\opens=msym8 %scaled\magstep1 % For Arbortxt on PC
  % For Arbortxt on PC, and VAX. Borde p199
\font\ssb=cmss10
\font\smsb=cmss8
\def\StretchRtArr#1{{\count255=0\loop\relbar\joinrel\advance\count255 by1
\ifnum\count255<#1\repeat\rightarrow}}
\def\StretchLtArr#1{\,{\leftarrow\!\!\count255=0\loop\relbar
\joinrel\advance\count255 by1\ifnum\count255<#1\repeat}}

\def\StretchLRtArr#1{\,{\leftarrow\!\!\count255=0\loop\relbar\joinrel\advance
\count255 by1\ifnum\count255<#1\repeat\rightarrow\,\,}}

\def\mbox#1{{\leavevmode\hbox{#1}}}

\def\hspace#1{{\phantom{\mbox#1}}}
\def\oZ{\mbox{\open\char90}}
\def\oC{\mbox{\open\char67}}

\def\ssi{\mbox{{\ssb\char105}}}
\def\ssj{\mbox{{\ssb\char106}}}

\def\ssA{\mbox{{\ssb\char65}}}
\def\smssA{\mbox{{\smsb\char65}}}
\def\ssB{\mbox{{\ssb\char66}}}
\def\smssB{\mbox{{\smsb\char66}}}

\def\ssE{\mbox{{\ssb\char69}}}
\def\dfone{\mbox{{\ssb\char49}}}
\def\sstwo{\mbox{{\ssb\char50}}}
\def\ssthree{\mbox{{\ssb\char51}}}
\def\ssfour{\mbox{{\ssb\char52}}}
\def\smssfour{\mbox{{\smsb\char52}}}
\def\ssfive{\mbox{{\ssb\char53}}}
\def\sssix{\mbox{{\ssb\char54}}}

 %in jyTeX
 %in jyTeX
 %in jyTeX
 %in jyTeX
 %in jyTeX
 %in jyTeX
 %in jyTeX
 %in jyTeX
 %in jyTeX
 %in jyTeX
 %in jyTeX
 %in jyTeX
% in jyTeX
% in jyTeX
% in jyTeX
% in jyTeX
% in jyTeX

\def\ga{\gamma}

\def\Ga{\Gamma}

\def\la{\lambda}

\def\om{\omega}

\def\si{\sigma}

\def\De{\Delta}

\def\caG{{\cal G}}

\def\caS{{\cal S}}

\def\zf{$\zeta$--function}

     % Newline

\def\frac#1/#2{\leavevmode\kern.1em
\raise.5ex\hbox{\the\scriptfont0 #1}\kern-.1em/\kern-.15em
\lower.25ex\hbox{\the\scriptfont0 #2}}
\def\sfrac#1/#2{\leavevmode\kern.1em
\raise.5ex\hbox{\the\scriptscriptfont0 #1}\kern-.1em/\kern-.15em
\lower.25ex\hbox{\the\scriptscriptfont0 #2}}

\def\gtorder{\mathrel{\raise.3ex\hbox{$>$}\mkern-14mu
             \lower0.6ex\hbox{$\sim$}}}
\def\ltorder{\mathrel{\raise.3ex\hbox{$<$}\mkern-14mu
             \lower0.6ex\hbox{$\sim$}}}

\def\semidirprod{\rlap{\ss C}\raise1pt\hbox{$\mkern.75mu\times$}}
\def\for{\lower6pt\hbox{$\Big|$}}
\def\fish{\kern-.25em{\phantom{abcde}\over \phantom{abcde}}\kern-.25em}

 %triple
%dot
 %double
%dot
 %double dot
%for small #1

\def\boxit#1{\vbox{\hrule\hbox{\vrule\kern3pt
        \vbox{\kern3pt#1\kern3pt}\kern3pt\vrule}\hrule}}
\def\dalemb#1#2{{\vbox{\hrule height .#2pt
        \hbox{\vrule width.#2pt height#1pt \kern#1pt \vrule
                width.#2pt} \hrule height.#2pt}}}

\def\ol{\overline}
        %double stroke
\def\frac#1#2{{{#1}\over{#2}}}
 %lower covariant deriv.
 %upper covariant deriv.
 %lower covariant deriv semicolon.
    %lower ordinary  deriv.
    %lower ordinary  deriv comma.

\def\noin{\noindent}
\def\Hom{\rm Hom}

      %Connection
    %Connection'

\def\etc{{\it etc. }}

\def\eg{{\it e.g.}}
\def\ie{{\it i.e. }}

\def\ket#1{\mid#1\rangle}
\def\kett#1{\mid#1\rangle\rangle}
 %gives average <#1>
 %gives thermal average <<#1>>
\def\br#1#2{\langle{#1}\mid{#2}\rangle}   %gives bracket <#1|#2>
   %gives comma bracket <#1,#2>
 %gives round bracket (#1,#2)
 %gives round bracket (#1,|#2)
 %gives big bracket <#1|#2>
  %gives
%matrix element <#1|#2|#3>
  %gives reduced matrix element
%<#1||#2||#3>

\def\wh{\widehat}
\def\wt{\widetilde}

\def\3j#1#2#3#4#5#6{\left\lgroup\matrix{#1&#2&#3\cr#4&#5&#6\cr}
\right\rgroup}

\def\man{{\cal M}}

\def\m?{\mgn{?}}
% KK's defs

\def\beq{\begin{eqnarray}}
\def\eeq{\end{eqnarray}}

%  *******************  Journal refs **********************

\def\aop#1#2#3{{\it Ann. Phys.} {\bf {#1}} ({#2}) #3}

\def\cmp#1#2#3{{\it Comm. Math. Phys.} {\bf {#1}} ({#2}) #3}
\def\cqg#1#2#3{{\it Class. Quant. Grav.} {\bf {#1}} ({#2}) #3}

\def\ijmp#1#2#3{{\it Int. J. Mod. Phys.} {\bf {#1}} ({#2}) #3}

\def\jmp#1#2#3{{\it J. Math. Phys.} {\bf {#1}} ({#2}) #3}
\def\jpa#1#2#3{{\it J. Phys.} {\bf A{#1}} ({#2}) #3}
\def\lnm#1#2#3{{\it Lect. Notes Math.} {\bf {#1}} ({#2}) #3}

\def\np#1#2#3{{\it Nucl. Phys.} {\bf B{#1}} ({#2}) #3}
\def\pl#1#2#3{{\it Phys. Lett.} {\bf {#1}} ({#2}) #3}

\def\prp#1#2#3{{\it Phys. Rep.} {\bf {#1}} ({#2}) #3}
\def\pr#1#2#3{{\it Phys. Rev.} {\bf {#1}} ({#2}) #3}
\def\prA#1#2#3{{\it Phys. Rev.} {\bf A{#1}} ({#2}) #3}

\def\prD#1#2#3{{\it Phys. Rev.} {\bf D{#1}} ({#2}) #3}
\def\prl#1#2#3{{\it Phys. Rev. Lett.} {\bf #1} ({#2}) #3}

\def\rmp#1#2#3{{\it Rev. Mod. Phys.} {\bf {#1}} ({#2}) #3}

\def\zfp#1#2#3{{\it Z. f. Phys.} {\bf {#1}} ({#2}) #3}

\def\cras#1#2#3{{\it Comptes Rend. Acad. Sci. (Paris)} {\bf{#1}} (#2) #3}
\def\prs#1#2#3{{\it Proc. Roy. Soc.} {\bf A{#1}} ({#2}) #3}
\def\pcps#1#2#3{{\it Proc. Camb. Phil. Soc.} {\bf{#1}} ({#2}) #3}
\def\mpcps#1#2#3{{\it Math. Proc. Camb. Phil. Soc.} {\bf{#1}} ({#2}) #3}

\def\amsh#1#2#3{{\it Abh. Math. Sem. Ham.} {\bf {#1}} ({#2}) #3}
\def\am#1#2#3{{\it Acta Mathematica} {\bf {#1}} ({#2}) #3}
\def\aim#1#2#3{{\it Adv. in Math.} {\bf {#1}} ({#2}) #3}
\def\ajm#1#2#3{{\it Am. J. Math.} {\bf {#1}} ({#2}) #3}

\def\aom#1#2#3{{\it Ann. of Math.} {\bf {#1}} ({#2}) #3}
\def\cjm#1#2#3{{\it Can. J. Math.} {\bf {#1}} ({#2}) #3}
\def\bams#1#2#3{{\it Bull.Am.Math.Soc.} {\bf {#1}} ({#2}) #3}

\def\cmh#1#2#3{{\it Comm. Math. Helv.} {\bf {#1}} ({#2}) #3}

\def\dmj#1#2#3{{\it Duke Math. J.} {\bf {#1}} ({#2}) #3}
\def\invm#1#2#3{{\it Invent. Math.} {\bf {#1}} ({#2}) #3}

\def\jdg#1#2#3{{\it J. Diff. Geom.} {\bf {#1}} ({#2}) #3}

\def\joa#1#2#3{{\it J. of Algebra} {\bf {#1}} ({#2}) #3}
\def\jram#1#2#3{{\it J. f. reine u. Angew. Math.} {\bf {#1}} ({#2}) #3}
\def\jims#1#2#3{{\it J. Indian. Math. Soc.} {\bf {#1}} ({#2}) #3}
\def\jlms#1#2#3{{\it J. Lond. Math. Soc.} {\bf {#1}} ({#2}) #3}
\def\jmpa#1#2#3{{\it J. Math. Pures. Appl.} {\bf {#1}} ({#2}) #3}
\def\ma#1#2#3{{\it Math. Ann.} {\bf {#1}} ({#2}) #3}

\def\mz#1#2#3{{\it Math. Zeit.} {\bf {#1}} ({#2}) #3}
\def\ojm#1#2#3{{\it Osaka J.Math.} {\bf {#1}} ({#2}) #3}

\def\pems#1#2#3{{\it Proc. Edin. Math. Soc.} {\bf {#1}} ({#2}) #3}

\def\plb#1#2#3{{\it Phys. Letts.} {\bf {B#1}} ({#2}) #3}
\def\pla#1#2#3{{\it Phys. Letts.} {\bf {A#1}} ({#2}) #3}
\def\plms#1#2#3{{\it Proc. Lond. Math. Soc.} {\bf {#1}} ({#2}) #3}
\def\pgma#1#2#3{{\it Proc. Glasgow Math. Ass.} {\bf {#1}} ({#2}) #3}
\def\qjm#1#2#3{{\it Quart. J. Math.} {\bf {#1}} ({#2}) #3}
\def\qjpam#1#2#3{{\it Quart. J. Pure and Appl. Math.} {\bf {#1}} ({#2}) #3}

\def\rmjm#1#2#3{{\it Rocky Mountain J. Math.} {\bf {#1}} ({#2}) #3}

\def\tams#1#2#3{{\it Trans.Am.Math.Soc.} {\bf {#1}} ({#2}) #3}

% *******************   Main text *********************
\begin{title}
\vglue 0.5truein
%\righttext {MUTP/96/23}
%\righttext{hep-th/96}
\vskip15truept
%\leftline{\today}
%\vskip 30truept
\centertext {\Bigfonts \bf Group theory aspects of spectral }
\vskip7truept \vskip10truept\centertext{\Bigfonts \bf  problems on
spherical factors}
 \vskip 20truept
\centertext{J.S.Dowker\footnote{dowker@man.ac.uk}} \vskip 7truept
\centertext{\it Theory Group,} \centertext{\it School of Physics and
Astronomy,} \centertext{\it The University of Manchester,}
\centertext{\it Manchester, England} \vskip 7truept \centertext{}

\vskip 7truept

\vskip40truept
\begin{narrow}
The Ray--Singer isospectral theorem (1971) is applied to a general
spectral function for Laplacians of twisted $p$--forms (say) on
homogeneous Clifford--Klein  factors of the three--sphere. The inducing
formulae necessary to express any spectral quantity for any twisting in
terms of those for cyclic subgroups of the tetrahedral, octahedral and
icosahedral deck groups  are detailed. Further, Artin's theorem allows
the McKay correspondence to be obtained.

The isospectral theorem is shown to yield a derivation of the Sunada
construction which is equivalent to the later one by Pesce.
\end{narrow}
\vskip 5truept
%\righttext {August 1996}
\vskip 60truept
%\righttext{Typeset in \jyTeX}
\vfil
\end{title}
\pagenum=0
\newpage

\section{\bf 1. Introduction.}

This paper is a product of my ongoing interest in explicit calculations
of spectral quantities, such as the Casimir energy and effective action,
on spherical factors as examples of manageable manifolds of non--trivial
topology and vaguely physical significance. It is always possible, of
course, to go for generality, and treat symmetric spaces, however I
prefer to concentrate on rather specific examples, in particular on
factors of the three--sphere, not only because the techniques are
commonly available but also because the factor possibilities are more
extensive and interesting.

The situation I wish to address here is a fairly common one and is the
same as that considered in [\pref{DandJ}] namely that of a (complex)
field, of some particular space--time character, belonging to a
representation of an internal summetry group, $G$, and defined on a
factor, S$^3/\Ga$ (to be specified later). Space--time could be $T\times
S^3/\Ga$ but it is the spectral problem on the spherical factor that I
propose to concentrate on and I henceforth ignore any time dependence.

In [\pref{DandJ}] we computed the Casimir energy, for various fields,
representations and factors. The calculations for $\Ga$ one of the binary
polyhedral groups, $T'$, $O'$, $Y'$, were performed individually. The
basic principle that I wish to investigate here is that all computations
can be reduced to those for cyclic groups.

I have implemented this earlier, [\pref{Dow11}], but only for untwisted
fields (and homogeneous factors). That it is possible in general follows
from Artin's theorem and an isospectral result contained in Ray and
Singer, [\pref{RandS}], applied to the analytic torsion, a specific
spectral quantity, but valid generally. I explain this in the next
section.
\section{\bf 2. The setup.}

The standard, generic setup, \eg\ [{\pref{RandS}], is a field, $\wt\phi$,
defined on the simply connected universal covering space, $\wt\man$,
satisfying the periodicity conditions (twisting),
  $$
  \wt\phi(x\ga)=\wt\phi(x)\,\rho(\ga)
  \eql{twist}
  $$
which projects down to a `multivalued field', $\phi$, on
$\man=\wt\man/\Ga$. The matrix, $\rho$, is a representation of the space
group, $\Ga$, in the internal group, $\caG$, \ie $\rho\in \Hom(\Ga,\caG)$
and $\rho(\ga\ga')=\rho(\ga)\rho(\ga')$. For concreteness I take $\caG$
to be U$(N)$ and $\phi$ in the fundamental representation so that $\rho$
is, initially, an $N\times N$ matrix, $||\rho_{ij}||$. Also, I will
choose $\phi$ to be a space-time scalar field, or, possibly, a $p$--form.
The set--up is a particular case of a more general situation, termed an
`automorphic' field theory in [\pref{BandD}]. In the present paper,
$\wt\phi$ is a section of a {\it flat} vector bundle which implies,
technically, that it is possible to choose frames such that the
heat--kernel, for example, on $\wt\man$ is proportional to the unit
matrix in internal (fibre) space.

A (Laplacian) spectral quantity, $\caS(\man;\rho)$, is defined to be a
function of the spectrum, $\{\la_n(\rho)\}$, of the de Rham Laplacian
\footnote{ The results of this paper are actually valid for any natural
operator.} on $\man$ for fields satisfying the twisting (\peq{twist}).
Examples might be the fully traced heat--kernel and \zf. For a flat
vector bundle these quantities involve only the character of the
twisting, $\rho$, and not the complete representation. To this fact can
be traced the computational tractability.

Ray and Singer, [\pref{RandS}], prove an `isospectral' theorem,
  $$
  \caS(\man_1;{\rm Ind}\,\rho)=\caS(\man_2;\rho)\,,
  \eql{isosp}
  $$
where $\man_1$ and $\man_2$ are both covered by $\wt\man$ and
$\Ga_2\subset\Ga_1$. They apply it to the analytic torsion, but it holds
in general. I give a quick proof in the next section, for completeness
and then apply this, rather simple, result to the spherical factors
mentioned above with $\Ga_1$ one of $T'$, $O'$ or $Y'$ and $\Ga_2$ a
cyclic subgroup (when $\man_2$ is a lens space and, therefore,
`simpler').

To explain the bookkeeping, I remark, again, that the representation
$\rho_1$ is $N$--dimensional and that an element of $\Hom(\Ga_1,\caG)$
can be specified by populating the $N\times N$ matrix with sufficient and
suitable irreps of $\Ga_1$, which are known. For example, there are five
non--trivial elements of ${\Hom}\big(Y',U(4)\big)$, corresponding to the
reps ${\ssfour}$, ${\ssfour_s}$,  ${\sstwo_s}\oplus{\sstwo_s}$,
${\sstwo'_s}\oplus{\sstwo'_s}$ and ${\sstwo_s}\oplus{\sstwo'_s}$ of $Y'$.
Therefore the elements corresponding to the individual irreps of $\Ga_1$
are building blocks from which any required $\Hom$ can be assembled.
These are the particular objects I seek.

If $\Ga_2$ is a cyclic subgroup, ${\rm Ind}\,\rho={\rm Ind}\,\om$ (where
$\om$ is a root of unity) and is (equivalent to) a direct sum of irreps
of $\Ga_1$. Putting this into (\peq{isosp}), the direct sum becomes an
algebraic one\footnote{ Using the assumed additivity,
$\caS(\man;\rho\oplus\rho')=\caS(\man;\rho)+\caS(\man;\rho')$, which
holds, \eg, for the \zf.} and so, turning it around, the spectral
quantity building block, $\caS(\man_1;\ssA)$ ($\ssA$ is an irrep of
$\Ga_1$), can be obtained as a linear combination of the
$\caS(\man_2;\om)$, the algebraic sufficiency of the cyclic quantities
being guaranteed by Artin's theorem. It is often possible to obtain these
latter in closed form.

It is my intention in this paper just to detail these linear
combinations. For simplicity, I consider only one--sided (homogeneous)
factors of the three--sphere.
\section{\bf 3. Isospectrality.}

The eigenvalues of the Laplacian on $\wt\man/\Ga$ are the same as those
on the covering manifold, $\wt\man$, only the degeneracies might differ
and the isospectrality, (\peq{isosp}), is really just a statement about
these.

As is well known in quantum mechanics, (Landau and Lifshitz,
[\pref{LandL}]), and many other areas, the eigenspace spanned by the
eigenvectors with a particular eigenvalue, $\la$, on $\wt\man$, forms the
carrier space of a rep of {\it any} symmetry group that $\wt\man$ might
possess, induced on the eigenvectors by the action of this group on
$\wt\man$.

Barring accidents, this rep is an irrep of the biggest symmetry group of
$\wt\man$ which, in the case of S$^3$, is O(4), although this plays no
role in the following. The dimension, $d_\la$, of this irrep,
$\wt\ssE_\la$, is the degeneracy of the corresponding energy level, still
referring to the Laplacian on $\wt\man$ and ignoring any internal degrees
of freedom.

For a chain of subgroups, $\Ga_1\supset\Ga_2\supset\ldots$, each such rep
is subduced from the previous one. (As in crystal field theory, I am
thinking of the groups, $\Ga_i$, as finite ones.) The dimension of {\it
all} these reps is $d_\la$,
  $$
  d_\la=\dim \wt\ssE_\la=\dim {\rm Sub}\,\wt\ssE_\la=
  \dim {\rm Sub}\,{\rm Sub}\wt\ssE_\la  =\ldots\,.
  $$

For ease, denote the rep of $\Ga_1$ by $\ssE_\la\equiv{\rm
Sub}\,\wt\ssE_\la$ and assume that, on $\man_2$, there is a field twisted
by (\peq{twist}) with $\rho$ now an irrep, $\ssB$. The rep Sub $\ssE_\la$
(of $\Ga_2$) decomposes according to all the possible irreps the
frequency of $\ssB$ being the dimension of the {\it twisted} eigenspace
on $\man_2$, \ie the twisted degeneracy on $\man_2$,
  $$
  d^{\man_2}_\la(\ssB)=\br{\ssB}{{\rm Sub}\,\ssE_\la\,}_{\Ga_2}\,.
  $$

The total degeneracy give the dimension relation,
  $$
  d_\la=\sum_i d^{\man_2}_\la(\ssB_i)\,,
  $$
summed over all irreps, $\ssB_i$, in $\ssE_\la$. There is a similar
formula for every subgroup.

If $\rho$ is not irreducible, the degeneracy is the intertwining number,
  $$
  d^{\man_2}_\la(\rho)=\br{\rho}{{\rm Sub}\,\ssE_\la\,}_{\Ga_2}\,,
  $$
by linearity of characters. The standard formula for the bracket is given
in the next section.

As noted by Ray and Singer, [\pref{RandS}], Frobenius reciprocity allows
one to write this as,
  $$
  \br{\rho}{{\rm Sub}\,\ssE_\la}_{\Ga_2}=
  \br{{\rm Ind}\,\rho}{\ssE_\la}_{\Ga_1}\,,
  $$
which is recognised as the degeneracy of the $\la$ eigenspace on $\man_1$
for fields twisted by Ind $\rho$ and so,
  $$
  d^{\man_2}_\la(\rho)=d^{\man_1}_\la({\rm Ind}\,\rho)\,.
  $$
This is the required result and leads to (\peq{isosp}).

This analysis is a purely representation--theoretic matter, like the
Sunada approach to isospectrality, [\pref{Sunada}], and its
generalisations by Pesce, [\pref{Pesce}], [\pref{GoandM}], and Sutton,
[\pref{Sutton}], and I will now give a derivation of this based on
({\peq{isosp}) which is essentially the same as Pesce's,
[\pref{Pesce2}].\footnote { The workers in this field seem to be unaware
of the relevance of the Ray--Singer theorem.}

The Sunada construction uses three (finite) groups $\Ga,\,\Ga_1$ and
$\Ga_2$, with $\Ga_1$ and $\Ga_2$ subgroups of $\Ga$, and all taken as
symmetry groups of some covering manifold $\wt\man$. Define, as above,
  $$
  \man=\wt\man/\Ga\,,\quad \man_1=\wt\man/\Ga_1\,,
  \quad \man_2=\wt\man/\Ga_2\,,
  $$
and apply (\peq{isosp}) to the two pairs, ($\man,\,\man_1$) and ($\man,\,
\man_2$) to give
  $$\eqalign{
  \caS(\man;{\rm Ind}\,\rho_1)&=\caS(\man_1;\rho_1)\cr
   \caS(\man;{\rm Ind}\,\rho_2)&=\caS(\man_2;\rho_2)
   }
  \eql{isosp2}
  $$
so that, if ${\rm Ind}\,\rho_1$ and ${\rm Ind}\,\rho_2$ are
$\Ga$--equivalent, one obtains the `isospectral' relation,
  $$
  \caS(\man_1;\rho_1)=\caS(\man_2;\rho_2)\,,
  \eql{isosp3}
  $$
generalising Sunada's original theorem, to which it reduces when $\rho_1$
and $\rho_2$ are the trivial reps making the equivalence condition
somewhat restrictive, but more significant in that the vector bundles are
`the same'.

\section{\bf 4. Inducing representations.}

There is a standard technique for computing representations of a group,
$G$, from those of a subgroup, $H$, which is especially easy if $H$ is
abelian. Given the character tables, routine algebra will produce the
answer.

Some basic, textbook facts and organising notation are necessary. I
denote the cyclic subgroups of the binary polyhedral groups generically
by $\oZ_q$. The non--trivial irreps of $\oZ_q$ are generated by the
$r$--powers of a primitive $q$--th root of unity where $1\le r\le q-1$.
The reps induced by these irreps are denoted initially by ${\rm
Ind}\,(\om_q^r)$, where $\om_q=e^{2\pi i/q}$. (It is not necessary to
consider subgroups, if any, of $\oZ_q$, independently, since inducing is
transitive.\footnote{ Mackey, [\pref{Mackey}], calls this `inducing in
stages'. Sometimes the phrase `inducing through' principle is used. })

The necessary values of the order, $q$, are contained in the
Threlfall--Coxeter presentation of the binary polyhedral groups, $\langle
l,m,n\rangle$
  $$
  \langle l,m,n\rangle: R^l=S^m=T^n=RST\,,\quad l=2,\,m=3,\,n=3,4,5\,.
  \eql{pres}
  $$
These imply that $(RST)^2=E={\rm id}$, \eg\ Coxeter and Moser,
[\pref{CandM}],\S6.5 so the orders, $q$, are $2l,2m$ and $2n$.

Because the $\oZ_q$ meet all the conjugacy classes of $\langle
l,m,n\rangle$ they are sufficient, by Artin's theorem, to induce all its
irreps.

The sufficiency of the irreps of the $\oZ_q$ can also be checked in the
following way. The number of irreps of $ \langle l,m,n\rangle$ is the
same as that of the conjugacy classes and these comprise the three,
`non--trivial' types,
  $$\eqalign{
  [R^j]&:\quad 1\le j\le l-1\,\cr
  [S^j]&:\quad 1\le j\le m-1\,\cr
  [T^j]&:\quad 1\le j\le n-1\,,\cr
  }
  \eql{nontrivcl}
  $$
(noting that $[T^j]=[T^{2n-j}]=[T^{-j}]$, \etc\footnote{ The tetrahedral
case, $\langle 2,3,3\rangle$, is more involved in that the two larger
class sets are {\it cross--linked}, \eg\ $[S]=[T^{-1}]$. I give some
details in  Appendix A, plus some other information.}) plus the two
trivial classes, $[E]$ and $[RST]=[\ol E]$, which each contain just a
single element. This counting is, self--evidently, exactly that of the
effectively independent irreps of the cyclic subgroups, generated by $R$,
$S$ and $T$. Expanding on this a little; the elements (classes), $T^j$
and $T^{2n-j}$ of the {\it subgroup}, $\approx \oZ_{2n}$, generated by
$T$, lie in the same class $[T^j]$ of $\langle l,m,n\rangle$. The
non-trivial irreps of $\oZ_{2n}$, $T\to\om_{2n}^r$ and
$T\to\om_{2n}^{2n-r}$ therefore give the same induced character, and so
are not distinct, in this regard. The counting of the distinct
non-trivial irreps is thus the same as that of the non-trivial classes,
(\peq{nontrivcl}). The trivial reps $T\to\om_{2n}^0\,(\equiv{\bf1})$ and
$T\to\om_{2n}^n\,(\equiv\ol{\bf1})$ are actually common to all cyclic
subgroups and contribute $1+1$ to the total number of distinct irreps of
the complete set, $\oZ_{2l}$, $\oZ_{2m}$ and $\oZ_{2n}$, which is
$1+1+l-1+m-1+n-1$, the usual value. This confirms Artin's theorem in this
case.

The construction of the induced reps is standard. Possibly it is easiest
to use Frobenius reciprocity which states that the number of times an
irrep, ${\ssA}$, of $G$ occurs in the rep induced from one, $\ssB$, of
$H$ is the same as the number of times $\ssB$ is contained is the rep of
$H$ subduced from $\ssA$. There is a standard formula for this frequency
as the scalar product (on $H$) of the corresponding characters
  $$
  n\big(\ssA,{\rm Ind}\,\ssB\big)=\br{\,\ssA}{{\rm Ind}\,\ssB}_G
  =\br{\,{\rm Sub}\,\ssA}{\ssB}_H={1\over|H|}
  \sum_h \ol\chi^{{\rm Sub}\,\smssA}(h)\,\chi^{\,\smssB}(h)\,.
  $$
Here, $H$ is $\oZ_q$ generated by $R,S,T$ in turn, and calculation
produces the following inductions.

\vglue .5truein

\noin For the tetrahedral case:

 $$
 \matrix{
 &&T\,\Leftrightarrow\,S^{-1}&&&&R\cr
 0\uparrow&=&{\dfone}+{\ssthree}&&&=&{\dfone}+{\dfone'}+{\dfone''}+{\ssthree}\cr
 1\uparrow&=&\!\sstwo_s''+\sstwo_s &&&=&\sstwo_s+\sstwo'_s+\sstwo''_s\cr
 2\uparrow&=&\,\,\dfone''+\ssthree &&&=&\,\,\ssthree+\ssthree\cr
 3\uparrow&=&\sstwo'_s+\sstwo''_s&&&&-\cr
 4\uparrow&=&\!\!\dfone'+\ssthree&&&&-\cr
 5\uparrow&=&\!\sstwo_s+\sstwo'_s&&&&-\,.\cr
 }
 $$

\vglue .5truein

\noin For the octahedral case:

 $$\matrix{
 &&T&&S&&R\cr
 0\uparrow &=&{\dfone}+\sstwo+{\ssthree}
 &=&{\dfone}+{\dfone'}+{\ssthree}+{\ssthree'}&=&
 {\dfone}+\sstwo+{\ssthree}+2\times{\ssthree'}\cr
 1\uparrow&=&\sstwo_s+\ssfour_s &=&\sstwo_s+\sstwo'_s+\ssfour_s
 &=&\sstwo_s+\sstwo'_s+2\times\ssfour_s\cr
 2\uparrow&=&\ssthree+\ssthree' &=&\,\,\sstwo+\ssthree+\ssthree'
 &=&\dfone'+\sstwo+2\times\ssthree+\ssthree'\cr
 3\uparrow&=&\sstwo'_s+\ssfour_s&=&\,\,\ssfour_s+\ssfour_s&&-\cr
 4\uparrow&=&\dfone'+\sstwo+\ssthree'&&-&&-\,.\cr
  }
$$
%\newpage\vfill

\noin For the icosahedral case:

$$\matrix{
 &&T&&S&&R\cr
 0\uparrow&=&{\dfone}+{\ssthree}+\ssthree'+\ssfive
 &=&{\dfone}+{\ssthree}+\ssthree'+2\!\times\ssfour+\ssfive
 &=&{\dfone}+{\ssthree}+\ssthree'+2\!\times\ssfour+3\!\times\ssfive\cr
 1\uparrow&=&\sstwo_s+\ssfour_s+\sssix_s
 &=&\sstwo_s+\sstwo'_s+\ssfour_s+2\times\sssix_s
 &=&\sstwo_s+\sstwo'_s+2\times\ssfour_s+3\times\sssix_s\cr
 2\uparrow&=&\ssthree'+\ssfour+\ssfive
 &=&\,\,\ssthree+\ssthree'+\ssfour+2\times\ssfive
 &=&2\times({\ssthree}+\ssthree'+\ssfour+\ssfive)\cr
 3\uparrow&=&\sstwo'_s+\ssfour_s+\sssix_s
  &=&\,\,2\times(\ssfour_s+\sssix_s)&&-\cr
 4\uparrow&=&\ssthree+\ssfour+\ssfive&&-&&-\cr
 5\uparrow&=&\sssix_s+\sssix_s&&-&&-\,.\cr
 }
 $$

The irreps of the groups, $T'$, $O'$, $Y'$, are labelled by their
dimension, distinguished by dashes and the spinor, double--valued ones
have a suffix `$s$' and the column entries cease when repetitions begin.
The notation now is that $r\!\uparrow$ refers to the rep induced by the
cyclic irrep generated by $\om^r$, $\om$ being a relevant primitive root
of unity; say $\om_{2n}$ for $T$, $\om_6$ for $S$ and $\om_4$ for $R$. To
specify the particular generator, if required, I write $3\!\uparrow \!T$
\etc often leaving the group implicit. $0\!\uparrow$ is sometimes
referred to as the {\it principal} induced rep. (See, \eg\ Lomont,
[\pref{Lomont}]) with special properties. For example, it contains the
trivial rep exactly once.

There are numerous checks of these results. For example one can induce to
$O'$ from $\oZ_6$ via $T'$. For example, $\big(2\uparrow
S\big)\uparrow\,=\,\big(\dfone''+\ssthree\big)\big|_{T'}\uparrow\,=
\sstwo+\ssthree+\ssthree'\big|_{O'}$ where, for convenience, I have used
the induction results from $T'$ to $O'$ listed in Stekolschchik,
[\pref{Stek}] p.178. Further, adding the complete columns gives the
regular representation and corresponds to inducing from the trivial rep
of $H=\{E\}$. As is well known, this is a consequence of Frobenius
reciprocity.

The previous counting shows that not all these relations are independent,
as can be confirmed visually.
\section{\bf 4. Spectral consequences}
The transition to spectral quantities, $\caS(\wt\man/\Ga;\rho)$, converts
the decompositions into algebraic equations which can be solved for the
$\caS(\wt\man/\Ga;\ssA)\equiv\caS(\ssA)$, where $\ssA$ is an irrep, in
terms of $\caS(\wt\man/\Ga;r\!\uparrow\! \ga)$. These, from
(\peq{isosp}), equal the lens space quantities
$\caS(\wt\man/\oZ_q;r)\equiv\caS(r;\ga)$, where $q$ is the order of the
generator $\ga,=R,S,T$. I recall that $r$ labels the twisting of the
$U(1)$ bundle on the lens space.

The spinor and non--spinor reps separate and in the case of
$\langle2,3,5\rangle$ elimination yields,
  $$
  \caS\!\left(\matrix{\sstwo_s\cr\sstwo'_s\cr\ssfour_s\cr\sssix_s}\right)=
  \left(\matrix{0&-1&-1/2&1\cr
                -1&0&-1/2&1\cr
                1&1&0&-1\cr
                0&0&1/2&0}\right)\,\caS\!
                \left(\matrix{1;T\cr 3;T
                \cr5; T\cr 1;S}
                \right)
                \eql{235s}
  $$
and
  $$
  \caS\!\left(\matrix{\dfone\cr\ssthree\cr\ssthree'\cr\ssfour\cr\ssfive}\right)=
  \left(\matrix{1&1&1&-1&-1/2\cr
                0&0&-1&0&1/2\cr
                0&-1&0&0&1/2\cr
                0&1&1&-1&0\cr
                0&0&0&1&-1/2}\right)\,\caS\!
                \left(\matrix{0;T\cr 2;T
                \cr 4;T\cr2; S\cr 2;R}
                \right)\,.
                \eql{235}
  $$

For $\langle2,3,4\rangle$,
  $$
  \caS\!\left(\matrix{\sstwo_s\cr\sstwo'_s\cr\ssfour_s}\right)=
  \left(\matrix{1&0&-1/2\cr
                -1&1&0\cr
                0&0&1/2}\right)\,\caS\!
                \left(\matrix{1;T\cr 1;S
                \cr 3;S}
                \right)
                \eql{234s}
  $$
and
  $$
  \caS\!\left(\matrix{\dfone\cr\dfone'\cr\sstwo\cr\ssthree\cr\ssthree'}\right)=
  \left(\matrix{1&1&1/2&-1&-1/2\cr
                0&0&1/2&-1&1/2\cr
                0&-1&0&1&0\cr
                0&0&-1/2&0&1/2\cr
                0&1&1/2&0&-1/2}\right)\,\caS\!
                \left(\matrix{0;T\cr 2;T
                \cr 4;T\cr2; S\cr2; R}
                \right)\,,
                \eql{234}
  $$
\vglue .25in \noin while for $\langle2,3,3\rangle$,
  $$
  \caS\!\left(\matrix{\sstwo_s\cr\sstwo'_s\cr\sstwo''_s}\right)=
  \left(\matrix{1&-1&1\cr
                1&1&-1\cr
                -1&1&1}\right)\,\caS\!
                \left(\matrix{1;T\cr 3;T
                \cr 5;T\cr}
                \right)
                \eql{233s}
  $$
and
  $$
  \caS\!\left(\matrix{\dfone\cr\dfone'\cr\dfone''\cr\ssthree\cr}\right)=
  \left(\matrix{1&0&0&-1/2\cr
                0&0&1&-1/2\cr
                0&1&0&-1/2\cr
                0&0&0&1/2}\right)\,\caS\!
                \left(\matrix{0;T\cr 2;T
                \cr 4;T\cr2; R}
                \right)\,.
                \eql{233}
  $$

As an example of a consistency check, one of many, the decomposition,
$1\!\uparrow\! R=\sstwo_s+\sstwo'_s+2\times\ssfour_s$ in the octahedral
case (not used in the derivation of (\peq{234s})) becomes
$\caS(R;1)=\caS(S;1)+\caS(S;3)/2$. This relates, as an illustration, the
analytic torsions on $\oZ_4$ and $\oZ_6$ lens spaces and, using Ray's
formula, numerically is $2=1\times\sqrt4$. (The additive $\caS$ is the
logarithm of the torsion.)
\section{\bf 5. Use of the $\oZ_2$ subgroup.}
While the three cyclic subgroups generated by $R,S$ and $T$ are
sufficient, additional use of the $\oZ_2$ subgroup generated by the
central element, $\ol E=RST$, provides a more symmetrical formulation.

Inducing gives,
  $$\matrix{
 &&T'&&O'&&Y'\cr
 \!\!\!\!\!\!\!0\uparrow \!\!\!\!\!\!\!\!&=&\!\!\!\!\!\!\!\!{\dfone}+\dfone'\!
 +\dfone''
 \!+3\!\!\times\!{\ssthree}\!\!\!\!
 &=&\!\!\!\!\!{\dfone}+{\dfone'}\!+2\times\sstwo+3\!\!\times\!({\ssthree}\!
 +\!{\ssthree'})\!\!\!\!&=&\!\!\!\!
 {\dfone}\!+3\!\times({\ssthree}\!+\!{\ssthree'})\!+\!4\!\!\times\!\ssfour
 +5\!\!\times\!{\ssfive'}\cr
 \!\!\!\!\!1\uparrow\!\!\!\!&=&\!\!\!\!2\!\!\times\!(\sstwo_s+\sstwo'_s+\sstwo''_s)
 \!\!\!\!&=&\!\!\!\!2\!\times(\sstwo_s+\sstwo'_s+2\!\!\times\ssfour_s)
 \!\!\!\!&=&\!\!\!\!2\!\!\times\!(\sstwo_s+\sstwo'_s+2\!\!\times\!\ssfour_s
 +3\!\!\times\!\sssix_s)\cr
  }
  \eql{z2ind}
$$
which illustrates nicely the regular rep result,
$0_{\{E\}}\!\!\uparrow\,=\,0_{\{E,\ol E\}}\!\!\uparrow+1_{\{E,\ol
E\}}\!\!\uparrow$, mentioned before.

Furthermore, the value of $\caS$ evaluated for the trivial bundle,
$\caS({\bf1})$, can be expressed purely in terms of {\it untwisted} lens
space values. For {\it all} factors of S$^3$, apart from lens spaces
themselves, it is easily established from the above listings that
  $$
   \caS({\bf1})={1\over2}\big(\caS(0;T)+\caS(0;S)+\caS(0;R)-\caS(0;RST)\big)\,.
   \eql{cycdec}
  $$

This relation was derived in [\pref{Dow11}] using a geometric, cyclic
decomposition of the traced (untwisted) heat--kernel (or, equivalently,
the \zf) on orbifold factors of the two--sphere which was obtained
earlier in [\pref{ChandD}]. The information used is, of course, contained
in the symmetry groups. An application to analytic torsion was made in
[\pref{DandCh2}] (see also Tsuchiya, [\pref{Tsuchiya2}]) and to Casimir
energies in [\pref{Dow11}].

In addition to (\peq{cycdec}), it is readily found that the same
combination evaluated for the first (spinor) twisting, yields the result,
valid for $T'$, $O'$ and $Y'$,
  $$
 \caS(\sstwo_s) =\caS(1;R)+\caS(1;S)+\caS(1;T)-\caS(1;RST)\,.
 \eql{cycdecsp}
  $$
(For the tetrahedral case, $\caS(1;S)$ equals $\caS(5;T)$.)

In similar vein, I find some other universal relations,\footnote{ For the
tetrahedral case, $\caS(\smssfour_s)$ is zero.}
  $$\eqalign{
  \caS(\ssthree)
  &=\caS({\bf1})+\caS(2;R)+\caS(2;S)+\caS(2;T)-\caS(2;RST)\cr
  \caS(\ssfour_s)&=\caS(\sstwo_s)+\caS(3;R)+\caS(3;S)+\caS(3;T)-\caS(3;RST)\cr
  }
  $$
and also, just for $Y'$,
  $$\eqalign{
  \caS(\ssthree')&=\caS({\bf1})+\caS(4;R)+\caS(4;S)+\caS(4;T)-\caS(4;RST)\cr
  \caS(\sssix_s)&=\caS(\ssfour_s)+\caS(5;R)+\caS(5;S)+\caS(5;T)-\caS(5;RST)\,,
  \cr
  }
  $$
with trivial equalities, $\caS(i+4;R)=\caS(i;R)=\caS(4-i;R)$,
$\caS(2;RST)=\caS(0;RST)$ \etc

\section{\bf 6. The McKay correspondence}
Making the cyclic twisting, $r$, correspond to $j$ in (\peq{nontrivcl}),
gives a two--to--one correspondence between the irreps of the three
cyclic subgroups\footnote{ I do not distinguish between the irreps
generated by different primitive roots of unity.} and the conjugacy
classes, and thence the irreps of $\langle l,m,n\rangle$. In fact one can
go further and link up with the McKay correspondence in the following
fashion.

Represent, in the usual cyclotomic way, the inducing cyclic irrep
generators, $\om^r$,  by points on three {\it distinct} unit circles,
best pictured as great circles on a two--sphere, intersecting at the
common, trivial rep points, ${\bf1}$ and ${\ol{\bf1}}$, represented by
the north and south poles. Then, for each circle, identify a semicircle
with its reflection under a $\oZ_2$ orbifold action with ${\bf1}$ and
${\ol{\bf1}}$ as fixed points. The three resulting semicircles give a
graph with these points as two three--nodes connected by three arcs with
$l-1$, $m-1$ and $n-1$ two--nodes each. This is a compactification of the
extended Dynkin diagram for ${\wt E}_{n+3}$ obtained by linking the two
shorter arms to the `affine' node, ${\bf1}$. This seems, to me, a more
symmetrical arrangement. The possibility of adding the identity element
node to the end of each branch of the Dynkin diagram is noted by
Rossmann, [\pref{Rossmann}]. Identifying these nodes goes a little
further and is in keeping with the geometrical interpretation where both
$\{E\}$ and $\{\ol E\}$ correspond to rotations through $2\pi$ and is the
reason I refer to $\{\ol E\}$ as a trivial class and to ${\ol{\bf1}}$ as
a trivial rep.

This is not the standard form of the Mckay correspondence, which usually
labels the nodes by the equivalence classes of the irreps of $\langle
l,m,n\rangle$, but Artin's theorem demonstrates they are effectively the
same. Furthermore, the construction of the previous paragraph obviously
applies with the inducing cyclic irreps replaced by the conjugacy classes
of $\langle l,m,n\rangle$, $[R^j]$ \etc This yields the dual of the
standard correspondence.

The class version of the McKay correspondence has been encountered before
by Ito and Reid, [\pref{IandR}], and discussed in more detail  from an
algebraic geometry perspective by Brylinski, [\pref{Brylinski}], whose
Thm.4.1, gives the rules for constructing a graph which turns out to be a
Dynkin diagram. The role of the quaternion representation is played by a
`special' class, corresponding to an end vertex of the graph. Whether two
vertices (classes) are graphically connected depends on relations between
representatives of these two classes and the special one.

Suter, [\pref{Suter}], contains suitably labelled Dynkin diagrams and
other useful information.

\section{\bf 7. Concluding remarks.}

I have given formulae, equns. (\peq{235s}) to (\peq{233}), that enable
any spectral quantity for a flat, twisted vector bundle over tetrahedral,
octahedral and icosahedral space to be found from the corresponding
quantity on lens spaces with various twistings. Applications will be
dealt with elsewhere. For example, the results of Cisneros--Molina,
[\pref{CM}], on the $\eta$--invariant of twisted Dirac operators can be
obtained in a more direct fashion.

\section{\bf Appendix A. The tetrahedral classes.}

The coupling between the arms of equal length of the Dynkin diagram of
$\langle 3,3,2\rangle$ is expressed by the class equalities,
$[S]=[T^{-1}]$ and $[S^2]=[T^{-2}]$ which can be shown by exhibiting the
conjugation. For example, $S^{-1}=U^{-1}TU$, where $U$ has to be a group
element. In fact $U=T^{-1}RT$ which is proved using the presentation
relations, (\peq{pres}). Directly,
  $$\eqalign{
  U^{-1}TU&=T^{-1}R^{-1}T\,T\,T^{-1}RT=T^{-2}S^{-1}TST^2=T^4S^5TST^2\cr&
  =T^4S^7T=T^4ST=T^3S^2=S^5=S^{-1}\,,
  }
  $$
where I have used the relations $R=ST$, $TST=S^2$ and $T^3=S^3$.

Ito and Reid, [\pref{IandR}], treat the conjugacy relations using a
different presentation.

Rossmann, [\pref{Rossmann}] Lemma 2.2, has also considered this coupling
using a standard quaternion representation of the generators given, \eg\
, in Coxeter, [\pref{Coxeter2}]. The conjugation (by $\ssi$) stated in
[\pref{Rossmann}] appears to be in error. It should be by $\ssj$ . The
change from $\ssi$ to $\ssj$ corresponds to the conjugation (rotation) by
$T$ in the above definition of $U$. See Coxeter, [\pref{Coxeter2}],p.75.

This cross--linking means that the picture leading to the (compactified)
Dynkin diagram has to be slightly amended for the $\langle 2,3,3\rangle$
case so that the $\oZ_2$ action (which is an inversion involution) now
identifies a semicircle of one circle with a semicircle of the other. The
upshot is that the Dynkin diagram consists of a complete circle for $S$
(or $T$) and a semicircle for $R$.

There is no cross--linking for the octahedral and icosahedral cases. For
example one can show that $T=U^{-1}T^{-1}U$ where $U=SRS^{-1}$ by a
similar manipulation as above. These conjugacy relations have a
geometrical significance.

\section{\bf Appendix B. Induced representations and isopectrality again}

There are many treatments of the notion of induced representations, which
goes back to Frobenius. Most use the coset decomposition of $G$. I set up
the left one,
  $$
  G=g_1H+g_2H+\ldots+g_nH=\bigcup_{i=1}^n g_iH\,,\quad n=|G|/|H|\,.
  \eql{cosdec2}
  $$
The $g_i$ can be taken as the representatives of the cosets. If a
particular set of representatives is chosen, every group element, $g$,
can be written uniquely as $g=g_ih$ for some $g_i$ and $h\in H$.

An unfussy way of proceeding is the following. Consider the basis
vectors, $\ket {\ssB,m}$ of a rep, $\ssB$, of the subgroup, $H$, and {\it
define} the {\it new} vectors, $\kett{\ssB,i,m}$ by the object,
  $$
  \kett{\ssB,i,m}=g_i\ket {\ssB,m}\,,\quad i=1,\ldots,n\,,\quad
  m=1,\ldots,d\,,
  \eql{newv}
  $$
the linear space of which I show to be closed under action by $G$.
Consider
  $$\eqalign{
  g\kett{\ssB,i,m}&=gg_i\ket {\ssB,m}=g_jh\ket {\ssB,m}\cr
   &=g_j\ket{\ssB,m'}\,D^{\smssB}_{m'm}(h)\cr
  &=\,\kett{\ssB,j,m'\,}\,D^{\smssB}_{m'm}(h)\,,
  }
  \eql{manip}
  $$
where the {\it left} coset decomposition has been used. Hence the
vectors, (\peq{newv}), form the basis for the carrier space of a
representation of $G$ {\it induced} from $\ssB$, of $H$, and denoted
$\ssB(H)\uparrow G$ or $\ssB\uparrow$, for short.

One can extract the representation matrices of $\ssB\uparrow$ from
(\peq{manip}),
  $$
  D^{\smssB\uparrow}_{jm',im}(g)=D^{\smssB}_{m'm}(h)
  $$
where $h=g_j^{-1}\,g\,g_i$.

In the above manipulation, (\peq{manip}), $g_i$ and $g$ are given, then
$g_j$ and $h$ are uniquely determined (given the set of representatives).
However, if we write
   $$
  D^{\smssB\uparrow}_{jm',im}(g)=D^{\smssB}_{m'm}(g_j^{-1}\,g\,g_i)\,,
  $$
where now $g_i$, $g$ {\it and} $g_j$ are given, we have to ensure that
$g_j^{-1}\,g\,g_i$ belongs to $H$. This can be achieved by including a
`Kronecker delta--function'
   $$
  D^{\smssB\uparrow}_{jm',im}(g)=\si_{ji}(g)\,
  D^{\smssB}_{m'm}(g_j^{-1}\,g\,g_i)\,,
  \eql{indrep2}
  $$
where
  $$
  \si_{ji}(g)=\cases{1\,,\quad{\rm if}\,\, g_j^{-1}\,g\,g_i\in H\cr0\,,{\rm
  \,\,\,\,\quad otherwise\,,}
  }
  $$
or one could write,
  $$
  D^{\smssB\uparrow}_{jm',im}(g)=\dot D^{\smssB}_{m'm}(g_j^{-1}\,g\,g_i)\,,
  $$
with the same import.

I now give an alternative approach to the isospectrality, (\peq{isosp}).
The space of $p$--forms on $\man_2$ twisted by a rep, $D^\smssB$ can be
identified with that of forms, $\phi$, taking values in $\oC^d$, on the
covering manifold, $\wt\man$, satisfying (\peq{twist}), (I have dropped
the tilde),
   $$
  \phi(x\ga_2)=\phi(x)\,D^\smssB(\ga_2)\,\,,\quad \ga_2\in\Ga_2\,.
  \eql{twist3}
  $$
Similarly, the space of $p$--forms, $\wh\phi$, on $\man_1$ twisted by
$\uparrow\!\ssB$ and taking values in $\oC^{nd}$, is equivalent to that
of forms on $\wt\man$ satisfying,
  $$
  \wh\phi(x\ga)= \wh\phi(x)\,D^{\smssB\uparrow}(\ga)\,\,,\quad \ga\in\Ga_1\,,
  \eql{twist4}
  $$
where $D^{\smssB\uparrow}$ is defined by (\peq{indrep2}).

To show the equivalence of these two spaces, one constructs a
one--to--one mapping between them, [\pref{RandS}].

For a form $\phi$ on $\wt\man$, with values in $\oC^d$ let $S\phi$ be the
form on $\wt\man$ with values in $\oC^{nd}$, defined by
  $$
  \big(S\phi\big)(x)=\sum_{\oplus i}\phi(x\ga_i)
  $$
where the $\ga_i$ are the representatives of the left cosets of $\Ga_2$
in $\Ga_1$.

The action of $\Ga_1$ on $S\phi$, is
  $$\eqalign{
  \big(S\phi\big)(x\ga )&=\sum_{\oplus i}\phi(x\ga\ga_i )\cr
  &=\sum_{\oplus i}\phi(x\ga_j\ga_j^{-1}\ga_2\ga_i )
  }
  $$
for {\it any} $\ga_j$. Now let $j$ range over $1\to n$ and sum over $j$.
If $\ga_j^{-1}\ga\ga_i$ belongs to $\Ga_2$, then one can apply the action
(\peq{twist3}). If $\ga_j^{-1}\ga\ga_i$ does not belong to $\Ga_2$, one
would want the result to be zero. This can be achieved by extending the
action to all of $\Ga_1$ by using the $\dot D$, which vanishes for all
these other $\ga_j\,$s. Then
  $$\eqalign{
  \big(S\phi\big)(x\ga )
  &=\sum_{\oplus i}\sum_j\phi(x\ga_j)\dot D(\ga_j^{-1}\ga\ga_i)\cr
  &=(S\phi)(x)\,D^{\smssB\uparrow}(\ga)\,,
  }
  $$
and $S\phi$ obeys (\peq{twist4}). The map $S$ is therefore into.To show
it is also onto, one needs the converse. A form taking values in
$\oC^{nd}$ has the general structure $\wh\phi=\sum_{\oplus i=1}^n\phi_i$
where each $\phi_i$ takes values in $\oC^d$. If it satisfies
(\peq{twist4}), then the component $\phi_1$ satisfies (\peq{twist3}) and,
further, $\wh\phi=S\phi_1$. Hence $S\phi$ is everything.

The projection $S$ commutes with the Laplacian since the latter commutes
with the action of $\Ga_1$ and so $S$ preserves eigenspaces.

\vskip5truept

 \noin{\bf References.} \vskip5truept
\begin{putreferences}
   \ref{Kohler}{K\"ohler,K.: Equivariant Reidemeister torsion on
   symmetric spaces. Math.Ann. {\bf 307}, 57-69 (1997)}
   \ref{Kohler2}{K\"ohler,K.: Equivariant analytic torsion on ${\bf P^nC}$.
   Math.Ann.{\bf 297}, 553-565 (1993) }
   \ref{Kohler3}{K\"ohler,K.: Holomorphic analytic torsion on Hermitian
   symmetric spaces. J.Reine Angew.Math. {\bf 460}, 93-116 (1995)}
   \ref{Zagierzf}{Zagier,D. {\it Zetafunktionen und Quadratische
   K\"orper}, (Springer--Verlag, Berlin, 1981).}
   \ref{Stek}{Stekholschkik,R. {\it Notes on Coxeter transformations and the McKay
   correspondence.} (Springer, Berlin, 2008).}
   \ref{Pesce}{Pesce,H. \cmh {71}{1996}{243}.}
   \ref{Pesce2}{Pesce,H. {\it Contemp. Math} {\bf 173} (1994) 231.}
   \ref{Sutton}{Sutton,C.J. {\it Equivariant isospectrality
   and isospectral deformations on spherical orbifolds}, ArXiv:math/0608567.}
   \ref{Sunada}{Sunada,T. \aom{121}{1985}{169}.}
   \ref{GoandM}{Gornet,R, and McGowan,J. {\it J.Comp. and Math.}
   {\bf 9} (2006) 270.}
   \ref{Suter}{Suter,R. {\it Manusc.Math.} {\bf 122} (2007) 1-21.}
   \ref{Lomont}{Lomont,J.S. {\it Applications of finite groups} (Academic
   Press, New York, 1959).}
   \ref{DandCh2}{Dowker,J.S. and Chang,Peter {\it Analytic torsion on
   spherical factors and tessellations}, arXiv:math.DG/0904.0744 .}
   \ref{Mackey}{Mackey,G. {\it Induced representations}
   (Benjamin, New York, 1968).}
   \ref{Koca}{Koca, {\it Turkish J.Physics}.}
   \ref{Brylinski}{Brylinski, J-L., {\it A correspondence dual to McKay's}
    ArXiv alg-geom/9612003.}
   \ref{Rossmann}{Rossman,W. {\it McKay's correspondence
   and characters of finite subgroups of\break SU(2)} {\it Progress in Math.}
      Birkhauser  (to appear) .}
   \ref{JandL}{James, G. and Liebeck, M. {\it Representations and
   characters of groups} (CUP, Cambridge, 2001).}
   \ref{IandR}{Ito,Y. and Reid,M. {\it The Mckay correspondence for finite
   subgroups of SL(3,C)} Higher dimensional varieties, (Trento 1994),
   221-240, (Berlin, de Gruyter 1996).}
   \ref{BandF}{Bauer,W. and Furutani, K. {\it J.Geom. and Phys.} {\bf
   58} (2008) 64.}
   \ref{Luck}{L\"uck,W. \jdg{37}{1993}{263}.}
   \ref{LandR}{Lott,J. and Rothenberg,M. \jdg{34}{1991}{431}.}
   \ref{DoandKi} {Dowker.J.S. and Kirsten, K. {\it Analysis and Appl.}
   {\bf 3} (2005) 45.}
   \ref{dowtess1}{Dowker,J.S. \cqg{23}{2006}{1}.}
   \ref{dowtess2}{Dowker,J.S. {\it J.Geom. and Phys.} {\bf 57} (2007) 1505.}
   \ref{MHS}{De Melo,T., Hartmann,L. and Spreafico,M. {\it Reidemeister
   Torsion and analytic torsion of discs}, ArXiv:0811.3196.}
   \ref{Vertman}{Vertman, B. {\it Analytic Torsion of a  bounded
   generalized cone}, ArXiv:0808.0449.}
   \ref{WandY} {Weng,L. and You,Y., {\it Int.J. of Math.}{\bf 7} (1996)
   109.}
   \ref{ScandT}{Schwartz, A.S. and Tyupkin,Yu.S. \np{242}{1984}{436}.}
   \ref{AAR}{Andrews, G.E., Askey,R. and Roy,R. {\it Special functions}
  (CUP, Cambridge, 1999).}
   \ref{Tsuchiya}{Tsuchiya, N.: R-torsion and analytic torsion for spherical
   Clifford-Klein manifolds.: J. Fac.Sci., Tokyo Univ. Sect.1 A, Mathematics
   {\bf 23}, 289-295 (1976).}
   \ref{Tsuchiya2}{Tsuchiya, N. J. Fac.Sci., Tokyo Univ. Sect.1 A, Mathematics
   {\bf 23}, 289-295 (1976).}
  \ref{Lerch}{Lerch,M. \am{11}{1887}{19}.}
  \ref{Lerch2}{Lerch,M. \am{29}{1905}{333}.}
  \ref{TandS}{Threlfall, W. and Seifert, H. \ma{104}{1930}{1}.}
  \ref{RandS}{Ray, D.B., and Singer, I. \aim{7}{1971}{145}.}
  \ref{RandS2}{Ray, D.B., and Singer, I. {\it Proc.Symp.Pure Math.}
  {\bf 23} (1973) 167.}
  \ref{Jensen}{Jensen,J.L.W.V. \aom{17}{1915-1916}{124}.}
  \ref{Rosenberg}{Rosenberg, S. {\it The Laplacian on a Riemannian Manifold}
  (CUP, Cambridge, 1997).}
  \ref{Nando2}{Nash, C. and O'Connor, D-J. {\it Int.J.Mod.Phys.}
  {\bf A10} (1995) 1779.}
  \ref{Fock}{Fock,V. \zfp{98}{1935}{145}.}
  \ref{Levy}{Levy,M. \prs {204}{1950}{145}.}
  \ref{Schwinger2}{Schwinger,J. \jmp{5}{1964}{1606}.}
  \ref{Muller}{M\"uller, \lnm{}{}{}.}
  \ref{VMK}{Varshalovich.}
  \ref{DandWo}{Dowker,J.S. and Wolski, A. \prA{46}{1992}{6417}.}
  \ref{Zeitlin1}{Zeitlin,V. {\it Physica D} {\bf 49} (1991).  }
  \ref{Zeitlin0}{Zeitlin,V. {\it Nonlinear World} Ed by
   V.Baryakhtar {\it et al},  Vol.I p.717,  (World Scientific, Singapore, 1989).}
  \ref{Zeitlin2}{Zeitlin,V. \prl{93}{2004}{264501}. }
  \ref{Zeitlin3}{Zeitlin,V. \pla{339}{2005}{316}. }
  \ref{Groenewold}{Groenewold, H.J. {\it Physica} {\bf 12} (1946) 405.}
  \ref{Cohen}{Cohen, L. \jmp{7}{1966}{781}.}
  \ref{AandW}{Argawal G.S. and Wolf, E. \prD{2}{1970}{2161,2187,2206}.}
  \ref{Jantzen}{Jantzen,R.T. \jmp{19}{1978}{1163}.}
  \ref{Moses2}{Moses,H.E. \aop{42}{1967}{343}.}
  \ref{Carmeli}{Carmeli,M. \jmp{9}{1968}{1987}.}
  \ref{SHS}{Siemans,M., Hancock,J. and Siminovitch,D. {\it Solid State
  Nuclear Magnetic Resonance} {\bf 31}(2007)35.}
 \ref{Dowk}{Dowker,J.S. \prD{28}{1983}{3013}.}
 \ref{Heine}{Heine, E. {\it Handbuch der Kugelfunctionen}
  (G.Reimer, Berlin. 1878, 1881).}
  \ref{Pockels}{Pockels, F. {\it \"Uber die Differentialgleichung $\De
  u+k^2u=0$} (Teubner, Leipzig. 1891).}
  \ref{Hamermesh}{Hamermesh, M., {\it Group Theory} (Addison--Wesley,
  Reading. 1962).}
  \ref{Racah}{Racah, G. {\it Group Theory and Spectroscopy}
  (Princeton Lecture Notes, 1951). }
  \ref{Gourdin}{Gourdin, M. {\it Basics of Lie Groups} (Editions
  Fronti\'eres, Gif sur Yvette. 1982.)}
  \ref{Clifford}{Clifford, W.K. \plms{2}{1866}{116}.}
  \ref{Story2}{Story, W.E. \plms{23}{1892}{265}.}
  \ref{Story}{Story, W.E. \ma{41}{1893}{469}.}
  \ref{Poole}{Poole, E.G.C. \plms{33}{1932}{435}.}
  \ref{Dickson}{Dickson, L.E. {\it Algebraic Invariants} (Wiley, N.Y.
  1915).}
  \ref{Dickson2}{Dickson, L.E. {\it Modern Algebraic Theories}
  (Sanborn and Co., Boston. 1926).}
  \ref{Hilbert2}{Hilbert, D. {\it Theory of algebraic invariants} (C.U.P.,
  Cambridge. 1993).}
  \ref{Olver}{Olver, P.J. {\it Classical Invariant Theory} (C.U.P., Cambridge.
  1999.)}
  \ref{AST}{A\v{s}erova, R.M., Smirnov, J.F. and Tolsto\v{i}, V.N. {\it
  Teoret. Mat. Fyz.} {\bf 8} (1971) 255.}
  \ref{AandS}{A\v{s}erova, R.M., Smirnov, J.F. \np{4}{1968}{399}.}
  \ref{Shapiro}{Shapiro, J. \jmp{6}{1965}{1680}.}
  \ref{Shapiro2}{Shapiro, J.Y. \jmp{14}{1973}{1262}.}
  \ref{NandS}{Noz, M.E. and Shapiro, J.Y. \np{51}{1973}{309}.}
  \ref{Cayley2}{Cayley, A. {\it Phil. Trans. Roy. Soc. Lond.}
  {\bf 144} (1854) 244.}
  \ref{Cayley3}{Cayley, A. {\it Phil. Trans. Roy. Soc. Lond.}
  {\bf 146} (1856) 101.}
  \ref{Wigner}{Wigner, E.P. {\it Gruppentheorie} (Vieweg, Braunschweig. 1931).}
  \ref{Sharp}{Sharp, R.T. \ajop{28}{1960}{116}.}
  \ref{Laporte}{Laporte, O. {\it Z. f. Naturf.} {\bf 3a} (1948) 447.}
  \ref{Lowdin}{L\"owdin, P-O. \rmp{36}{1964}{966}.}
  \ref{Ansari}{Ansari, S.M.R. {\it Fort. d. Phys.} {\bf 15} (1967) 707.}
  \ref{SSJR}{Samal, P.K., Saha, R., Jain, P. and Ralston, J.P. {\it
  Testing Isotropy of Cosmic Microwave Background Radiation},
  astro-ph/0708.2816.}
  \ref{Lachieze}{Lachi\'eze-Rey, M. {\it Harmonic projection and
  multipole Vectors}. astro- \break ph/0409081.}
  \ref{CHS}{Copi, C.J., Huterer, D. and Starkman, G.D.
  \prD{70}{2003}{043515}.}
  \ref{Jaric}{Jari\'c, J.P. {\it Int. J. Eng. Sci.} {\bf 41} (2003) 2123.}
  \ref{RandD}{Roche, J.A. and Dowker, J.S. \jpa{1}{1968}{527}.}
  \ref{KandW}{Katz, G. and Weeks, J.R. \prD{70}{2004}{063527}.}
  \ref{Waerden}{van der Waerden, B.L. {\it Die Gruppen-theoretische
  Methode in der Quantenmechanik} (Springer, Berlin. 1932).}
  \ref{EMOT}{Erdelyi, A., Magnus, W., Oberhettinger, F. and Tricomi, F.G. {
  \it Higher Transcendental Functions} Vol.1 (McGraw-Hill, N.Y. 1953).}
  \ref{Dowzilch}{Dowker, J.S. {\it Proc. Phys. Soc.} {\bf 91} (1967) 28.}
  \ref{DandD}{Dowker, J.S. and Dowker, Y.P. {\it Proc. Phys. Soc.}
  {\bf 87} (1966) 65.}
  \ref{DandD2}{Dowker, J.S. and Dowker, Y.P. \prs{}{}{}.}
  \ref{Dowk3}{Dowker,J.S. \cqg{7}{1990}{1241}.}
  \ref{Dowk5}{Dowker,J.S. \cqg{7}{1990}{2353}.}
  \ref{CoandH}{Courant, R. and Hilbert, D. {\it Methoden der
  Mathematischen Physik} vol.1 \break (Springer, Berlin. 1931).}
  \ref{Applequist}{Applequist, J. \jpa{22}{1989}{4303}.}
  \ref{Torruella}{Torruella, \jmp{16}{1975}{1637}.}
  \ref{Weinberg}{Weinberg, S.W. \pr{133}{1964}{B1318}.}
  \ref{Meyerw}{Meyer, W.F. {\it Apolarit\"at und rationale Curven}
  (Fues, T\"ubingen. 1883.) }
  \ref{Ostrowski}{Ostrowski, A. {\it Jahrsb. Deutsch. Math. Verein.} {\bf
  33} (1923) 245.}
  \ref{Kramers}{Kramers, H.A. {\it Grundlagen der Quantenmechanik}, (Akad.
  Verlag., Leipzig, 1938).}
  \ref{ZandZ}{Zou, W.-N. and Zheng, Q.-S. \prs{459}{2003}{527}.}
  \ref{Weeks1}{Weeks, J.R. {\it Maxwell's multipole vectors
  and the CMB}.  astro-ph/0412231.}
  \ref{Corson}{Corson, E.M. {\it Tensors, Spinors and Relativistic Wave
  Equations} (Blackie, London. 1950).}
  \ref{Rosanes}{Rosanes, J. \jram{76}{1873}{312}.}
  \ref{Salmon}{Salmon, G. {\it Lessons Introductory to the Modern Higher
  Algebra} 3rd. edn. \break (Hodges,  Dublin. 1876.)}
  \ref{Milnew}{Milne, W.P. {\it Homogeneous Coordinates} (Arnold. London. 1910).}
  \ref{Niven}{Niven, W.D. {\it Phil. Trans. Roy. Soc.} {\bf 170} (1879) 393.}
  \ref{Scott}{Scott, C.A. {\it An Introductory Account of
  Certain Modern Ideas and Methods in Plane Analytical Geometry,}
  (MacMillan, N.Y. 1896).}
  \ref{Bargmann}{Bargmann, V. \rmp{34}{1962}{300}.}
  \ref{Maxwell}{Maxwell, J.C. {\it A Treatise on Electricity and
  Magnetism} 2nd. edn. (Clarendon Press, Oxford. 1882).}
  \ref{BandL}{Biedenharn, L.C. and Louck, J.D. {\it Angular Momentum in Quantum Physics}
  (Addison-Wesley, Reading. 1981).}
  \ref{Weylqm}{Weyl, H. {\it The Theory of Groups and Quantum Mechanics}
  (Methuen, London. 1931).}
  \ref{Robson}{Robson, A. {\it An Introduction to Analytical Geometry} Vol I
  (C.U.P., Cambridge. 1940.)}
  \ref{Sommerville}{Sommerville, D.M.Y. {\it Analytical Conics} 3rd. edn.
   (Bell. London. 1933).}
  \ref{Coolidge}{Coolidge, J.L. {\it A Treatise on Algebraic Plane Curves}
  (Clarendon Press, Oxford. 1931).}
  \ref{SandK}{Semple, G. and Kneebone. G.T. {\it Algebraic Projective
  Geometry} (Clarendon Press, Oxford. 1952).}
  \ref{AandC}{Abdesselam A., and Chipalkatti, J. {\it The Higher
  Transvectants are redundant}, arXiv:0801.1533 [math.AG] 2008.}
  \ref{Elliott}{Elliott, E.B. {\it The Algebra of Quantics} 2nd edn.
  (Clarendon Press, Oxford. 1913).}
  \ref{Elliott2}{Elliott, E.B. \qjpam{48}{1917}{372}.}
  \ref{Howe}{Howe, R. \tams{313}{1989}{539}.}
  \ref{Clebsch}{Clebsch, A. \jram{60}{1862}{343}.}
  \ref{Prasad}{Prasad, G. \ma{72}{1912}{136}.}
  \ref{Dougall}{Dougall, J. \pems{32}{1913}{30}.}
  \ref{Penrose}{Penrose, R. \aop{10}{1960}{171}.}
  \ref{Penrose2}{Penrose, R. \prs{273}{1965}{171}.}
  \ref{Burnside}{Burnside, W.S. \qjm{10}{1870}{211}. }
  \ref{Lindemann}{Lindemann, F. \ma{23} {1884}{111}.}
  \ref{Backus}{Backus, G. {\it Rev. Geophys. Space Phys.} {\bf 8} (1970) 633.}
  \ref{Baerheim}{Baerheim, R. {\it Q.J. Mech. appl. Math.} {\bf 51} (1998) 73.}
  \ref{Lense}{Lense, J. {\it Kugelfunktionen} (Akad.Verlag, Leipzig. 1950).}
  \ref{Littlewood}{Littlewood, D.E. \plms{50}{1948}{349}.}
  \ref{Fierz}{Fierz, M. {\it Helv. Phys. Acta} {\bf 12} (1938) 3.}
  \ref{Williams}{Williams, D.N. {\it Lectures in Theoretical Physics} Vol. VII,
  (Univ.Colorado Press, Boulder. 1965).}
  \ref{Dennis}{Dennis, M. \jpa{37}{2004}{9487}.}
  \ref{Pirani}{Pirani, F. {\it Brandeis Lecture Notes on
  General Relativity,} edited by S. Deser and K. Ford. (Brandeis, Mass. 1964).}
  \ref{Sturm}{Sturm, R. \jram{86}{1878}{116}.}
  \ref{Schlesinger}{Schlesinger, O. \ma{22}{1883}{521}.}
  \ref{Askwith}{Askwith, E.H. {\it Analytical Geometry of the Conic
  Sections} (A.\&C. Black, London. 1908).}
  \ref{Todd}{Todd, J.A. {\it Projective and Analytical Geometry}.
  (Pitman, London. 1946).}
  \ref{Glenn}{Glenn. O.E. {\it Theory of Invariants} (Ginn \& Co, N.Y. 1915).}
  \ref{DowkandG}{Dowker, J.S. and Goldstone, M. \prs{303}{1968}{381}.}
  \ref{Turnbull}{Turnbull, H.A. {\it The Theory of Determinants,
  Matrices and Invariants} 3rd. edn. (Dover, N.Y. 1960).}
  \ref{MacMillan}{MacMillan, W.D. {\it The Theory of the Potential}
  (McGraw-Hill, N.Y. 1930).}
   \ref{Hobson}{Hobson, E.W. {\it The Theory of Spherical and Ellipsoidal Harmonics}
   C.U.P., Cambridge. 1931).}
  \ref{Hobson1}{Hobson, E.W. \plms {24}{1892}{55}.}
  \ref{GandY}{Grace, J.H. and Young, A. {\it The Algebra of Invariants}
  (C.U.P., Cambridge, 1903).}
  \ref{FandR}{Fano, U. and Racah, G. {\it Irreducible Tensorial Sets}
  (Academic Press, N.Y. 1959).}
  \ref{TandT}{Thomson, W. and Tait, P.G. {\it Treatise on Natural Philosophy}
  (Clarendon Press, Oxford. 1867).}
  \ref{Brinkman}{Brinkman, H.C. {\it Applications of spinor invariants in
atomic physics}, North Holland, Amsterdam 1956.}
  \ref{Kramers1}{Kramers, H.A. {\it Proc. Roy. Soc. Amst.} {\bf 33} (1930) 953.}
  \ref{DandP2}{Dowker,J.S. and Pettengill,D.F. \jpa{7}{1974}{1527}}
  \ref{Dowk1}{Dowker,J.S. \jpa{}{}{45}.}
  \ref{Dowk2}{Dowker,J.S. \aop{71}{1972}{577}}
  \ref{DandA}{Dowker,J.S. and Apps, J.S. \cqg{15}{1998}{1121}.}
  \ref{Weil}{Weil,A., {\it Elliptic functions according to Eisenstein
  and Kronecker}, Springer, Berlin, 1976.}
  \ref{Ling}{Ling,C-H. {\it SIAM J.Math.Anal.} {\bf5} (1974) 551.}
  \ref{Ling2}{Ling,C-H. {\it J.Math.Anal.Appl.}(1988).}
 \ref{BMO}{Brevik,I., Milton,K.A. and Odintsov, S.D. \aop{302}{2002}{120}.}
 \ref{KandL}{Kutasov,D. and Larsen,F. {\it JHEP} 0101 (2001) 1.}
 \ref{KPS}{Klemm,D., Petkou,A.C. and Siopsis {\it Entropy
 bounds, monoticity properties and scaling in CFT's}. hep-th/0101076.}
 \ref{DandC}{Dowker,J.S. and Critchley,R. \prD{15}{1976}{1484}.}
 \ref{AandD}{Al'taie, M.B. and Dowker, J.S. \prD{18}{1978}{3557}.}
 \ref{Dow1}{Dowker,J.S. \prD{37}{1988}{558}.}
 \ref{Dow30}{Dowker,J.S. \prD{28}{1983}{3013}.}
 \ref{DandK}{Dowker,J.S. and Kennedy,G. \jpa{}{1978}{}.}
 \ref{Dow2}{Dowker,J.S. \cqg{1}{1984}{359}.}
 \ref{DandKi}{Dowker,J.S. and Kirsten, K. {\it Comm. in Anal. and Geom.
 }{\bf7} (1999) 641.}
 \ref{DandKe}{Dowker,J.S. and Kennedy,G.\jpa{11}{1978}{895}.}
 \ref{Gibbons}{Gibbons,G.W. \pl{60A}{1977}{385}.}
 \ref{Cardy}{Cardy,J.L. \np{366}{1991}{403}.}
 \ref{ChandD}{Chang,P. and Dowker,J.S. \np{395}{1993}{407}.}
 \ref{DandC2}{Dowker,J.S. and Critchley,R. \prD{13}{1976}{224}.}
 \ref{Camporesi}{Camporesi,R. \prp{196}{1990}{1}.}
 \ref{BandM}{Brown,L.S. and Maclay,G.J. \pr{184}{1969}{1272}.}
 \ref{CandD}{Candelas,P. and Dowker,J.S. \prD{19}{1979}{2902}.}
 \ref{Unwin1}{Unwin,S.D. Thesis. University of Manchester. 1979.}
 \ref{Unwin2}{Unwin,S.D. \jpa{13}{1980}{313}.}
 \ref{DandB}{Dowker,J.S.and Banach,R. \jpa{11}{1978}{2255}.}
 \ref{Obhukov}{Obhukov,Yu.N. \pl{109B}{1982}{195}.}
 \ref{Kennedy}{Kennedy,G. \prD{23}{1981}{2884}.}
 \ref{CandT}{Copeland,E. and Toms,D.J. \np {255}{1985}{201}.}
 \ref{ELV}{Elizalde,E., Lygren, M. and Vassilevich,
 D.V. \jmp {37}{1996}{3105}.}
 \ref{Malurkar}{Malurkar,S.L. {\it J.Ind.Math.Soc} {\bf16} (1925/26) 130.}
 \ref{Glaisher}{Glaisher,J.W.L. {\it Messenger of Math.} {\bf18}
(1889) 1.} \ref{Anderson}{Anderson,A. \prD{37}{1988}{536}.}
 \ref{CandA}{Cappelli,A. and D'Appollonio,\pl{487B}{2000}{87}.}
 \ref{Wot}{Wotzasek,C. \jpa{23}{1990}{1627}.}
 \ref{RandT}{Ravndal,F. and Tollesen,D. \prD{40}{1989}{4191}.}
 \ref{SandT}{Santos,F.C. and Tort,A.C. \pl{482B}{2000}{323}.}
 \ref{FandO}{Fukushima,K. and Ohta,K. {\it Physica} {\bf A299} (2001) 455.}
 \ref{GandP}{Gibbons,G.W. and Perry,M. \prs{358}{1978}{467}.}
 \ref{Dow4}{Dowker,J.S..}
  \ref{Rad}{Rademacher,H. {\it Topics in analytic number theory,}
Springer-Verlag,  Berlin,1973.}
  \ref{Halphen}{Halphen,G.-H. {\it Trait\'e des Fonctions Elliptiques},
  Vol 1, Gauthier-Villars, Paris, 1886.}
  \ref{CandW}{Cahn,R.S. and Wolf,J.A. {\it Comm.Mat.Helv.} {\bf 51}
  (1976) 1.}
  \ref{Berndt}{Berndt,B.C. \rmjm{7}{1977}{147}.}
  \ref{Hurwitz}{Hurwitz,A. \ma{18}{1881}{528}.}
  \ref{Hurwitz2}{Hurwitz,A. {\it Mathematische Werke} Vol.I. Basel,
  Birkhauser, 1932.}
  \ref{Berndt2}{Berndt,B.C. \jram{303/304}{1978}{332}.}
  \ref{RandA}{Rao,M.B. and Ayyar,M.V. \jims{15}{1923/24}{150}.}
  \ref{Hardy}{Hardy,G.H. \jlms{3}{1928}{238}.}
  \ref{TandM}{Tannery,J. and Molk,J. {\it Fonctions Elliptiques},
   Gauthier-Villars, Paris, 1893--1902.}
  \ref{schwarz}{Schwarz,H.-A. {\it Formeln und
  Lehrs\"atzen zum Gebrauche..},Springer 1893.(The first edition was 1885.)
  The French translation by Henri Pad\'e is {\it Formules et Propositions
  pour L'Emploi...},Gauthier-Villars, Paris, 1894}
  \ref{Hancock}{Hancock,H. {\it Theory of elliptic functions}, Vol I.
   Wiley, New York 1910.}
  \ref{watson}{Watson,G.N. \jlms{3}{1928}{216}.}
  \ref{MandO}{Magnus,W. and Oberhettinger,F. {\it Formeln und S\"atze},
  Springer-Verlag, Berlin 1948.}
  \ref{Klein}{Klein,F. {\it Lectures on the Icosohedron}
  (Methuen, London, 1913).}
  \ref{AandL}{Appell,P. and Lacour,E. {\it Fonctions Elliptiques},
  Gauthier-Villars,
  Paris, 1897.}
  \ref{HandC}{Hurwitz,A. and Courant,C. {\it Allgemeine Funktionentheorie},
  Springer,
  Berlin, 1922.}
  \ref{WandW}{Whittaker,E.T. and Watson,G.N. {\it Modern analysis},
  Cambridge 1927.}
  \ref{SandC}{Selberg,A. and Chowla,S. \jram{227}{1967}{86}. }
  \ref{zucker}{Zucker,I.J. {\it Math.Proc.Camb.Phil.Soc} {\bf 82 }(1977)
  111.}
  \ref{glasser}{Glasser,M.L. {\it Maths.of Comp.} {\bf 25} (1971) 533.}
  \ref{GandW}{Glasser, M.L. and Wood,V.E. {\it Maths of Comp.} {\bf 25}
  (1971)
  535.}
  \ref{greenhill}{Greenhill,A,G. {\it The Applications of Elliptic
  Functions}, MacMillan, London, 1892.}
  \ref{Weierstrass}{Weierstrass,K. {\it J.f.Mathematik (Crelle)}
{\bf 52} (1856) 346.}
  \ref{Weierstrass2}{Weierstrass,K. {\it Mathematische Werke} Vol.I,p.1,
  Mayer u. M\"uller, Berlin, 1894.}
  \ref{Fricke}{Fricke,R. {\it Die Elliptische Funktionen und Ihre Anwendungen},
    Teubner, Leipzig. 1915, 1922.}
  \ref{Konig}{K\"onigsberger,L. {\it Vorlesungen \"uber die Theorie der
 Elliptischen Funktionen},  \break Teubner, Leipzig, 1874.}
  \ref{Milne}{Milne,S.C. {\it The Ramanujan Journal} {\bf 6} (2002) 7-149.}
  \ref{Schlomilch}{Schl\"omilch,O. {\it Ber. Verh. K. Sachs. Gesell. Wiss.
  Leipzig}  {\bf 29} (1877) 101-105; {\it Compendium der h\"oheren
  Analysis}, Bd.II, 3rd Edn, Vieweg, Brunswick, 1878.}
  \ref{BandB}{Briot,C. and Bouquet,C. {\it Th\`eorie des Fonctions
  Elliptiques}, Gauthier-Villars, Paris, 1875.}
  \ref{Dumont}{Dumont,D. \aim {41}{1981}{1}.}
  \ref{Andre}{Andr\'e,D. {\it Ann.\'Ecole Normale Superior} {\bf 6} (1877)
  265;
  {\it J.Math.Pures et Appl.} {\bf 5} (1878) 31.}
  \ref{Raman}{Ramanujan,S. {\it Trans.Camb.Phil.Soc.} {\bf 22} (1916) 159;
 {\it Collected Papers}, Cambridge, 1927}
  \ref{Weber}{Weber,H.M. {\it Lehrbuch der Algebra} Bd.III, Vieweg,
  Brunswick 190  3.}
  \ref{Weber2}{Weber,H.M. {\it Elliptische Funktionen und algebraische
  Zahlen},
  Vieweg, Brunswick 1891.}
  \ref{ZandR}{Zucker,I.J. and Robertson,M.M.
  {\it Math.Proc.Camb.Phil.Soc} {\bf 95 }(1984) 5.}
  \ref{JandZ1}{Joyce,G.S. and Zucker,I.J.
  {\it Math.Proc.Camb.Phil.Soc} {\bf 109 }(1991) 257.}
  \ref{JandZ2}{Zucker,I.J. and Joyce.G.S.
  {\it Math.Proc.Camb.Phil.Soc} {\bf 131 }(2001) 309.}
  \ref{zucker2}{Zucker,I.J. {\it SIAM J.Math.Anal.} {\bf 10} (1979) 192,}
  \ref{BandZ}{Borwein,J.M. and Zucker,I.J. {\it IMA J.Math.Anal.} {\bf 12}
  (1992) 519.}
  \ref{Cox}{Cox,D.A. {\it Primes of the form $x^2+n\,y^2$}, Wiley,
  New York, 1989.}
  \ref{BandCh}{Berndt,B.C. and Chan,H.H. {\it Mathematika} {\bf42} (1995)
  278.}
  \ref{EandT}{Elizalde,R. and Tort.hep-th/}
  \ref{KandS}{Kiyek,K. and Schmidt,H. {\it Arch.Math.} {\bf 18} (1967) 438.}
  \ref{Oshima}{Oshima,K. \prD{46}{1992}{4765}.}
  \ref{greenhill2}{Greenhill,A.G. \plms{19} {1888} {301}.}
  \ref{Russell}{Russell,R. \plms{19} {1888} {91}.}
  \ref{BandB}{Borwein,J.M. and Borwein,P.B. {\it Pi and the AGM}, Wiley,
  New York, 1998.}
  \ref{Resnikoff}{Resnikoff,H.L. \tams{124}{1966}{334}.}
  \ref{vandp}{Van der Pol, B. {\it Indag.Math.} {\bf18} (1951) 261,272.}
  \ref{Rankin}{Rankin,R.A. {\it Modular forms} C.U.P. Cambridge}
  \ref{Rankin2}{Rankin,R.A. {\it Proc. Roy.Soc. Edin.} {\bf76 A} (1976) 107.}
  \ref{Skoruppa}{Skoruppa,N-P. {\it J.of Number Th.} {\bf43} (1993) 68 .}
  \ref{Down}{Dowker.J.S. \np {104}{2002}{153}.}
  \ref{Eichler}{Eichler,M. \mz {67}{1957}{267}.}
  \ref{Zagier}{Zagier,D. \invm{104}{1991}{449}.}
  \ref{Lang}{Lang,S. {\it Modular Forms}, Springer, Berlin, 1976.}
  \ref{Kosh}{Koshliakov,N.S. {\it Mess.of Math.} {\bf 58} (1928) 1.}
  \ref{BandH}{Bodendiek, R. and Halbritter,U. \amsh{38}{1972}{147}.}
  \ref{Smart}{Smart,L.R., \pgma{14}{1973}{1}.}
  \ref{Grosswald}{Grosswald,E. {\it Acta. Arith.} {\bf 21} (1972) 25.}
  \ref{Kata}{Katayama,K. {\it Acta Arith.} {\bf 22} (1973) 149.}
  \ref{Ogg}{Ogg,A. {\it Modular forms and Dirichlet series} (Benjamin,
  New York,
   1969).}
  \ref{Bol}{Bol,G. \amsh{16}{1949}{1}.}
  \ref{Epstein}{Epstein,P. \ma{56}{1903}{615}.}
  \ref{Petersson}{Petersson.}
  \ref{Serre}{Serre,J-P. {\it A Course in Arithmetic}, Springer,
  New York, 1973.}
  \ref{Schoenberg}{Schoenberg,B., {\it Elliptic Modular Functions},
  Springer, Berlin, 1974.}
  \ref{Apostol}{Apostol,T.M. \dmj {17}{1950}{147}.}
  \ref{Ogg2}{Ogg,A. {\it Lecture Notes in Math.} {\bf 320} (1973) 1.}
  \ref{Knopp}{Knopp,M.I. \dmj {45}{1978}{47}.}
  \ref{Knopp2}{Knopp,M.I. \invm {}{1994}{361}.}
  \ref{LandZ}{Lewis,J. and Zagier,D. \aom{153}{2001}{191}.}
  \ref{DandK1}{Dowker,J.S. and Kirsten,K. {\it Elliptic functions and
  temperature inversion symmetry on spheres} hep-th/.}
  \ref{HandK}{Husseini and Knopp.}
  \ref{Kober}{Kober,H. \mz{39}{1934-5}{609}.}
  \ref{HandL}{Hardy,G.H. and Littlewood, \am{41}{1917}{119}.}
  \ref{Watson}{Watson,G.N. \qjm{2}{1931}{300}.}
  \ref{SandC2}{Chowla,S. and Selberg,A. {\it Proc.Nat.Acad.} {\bf 35}
  (1949) 371.}
  \ref{Landau}{Landau, E. {\it Lehre von der Verteilung der Primzahlen},
  (Teubner, Leipzig, 1909).}
  \ref{Berndt4}{Berndt,B.C. \tams {146}{1969}{323}.}
  \ref{Berndt3}{Berndt,B.C. \tams {}{}{}.}
  \ref{Bochner}{Bochner,S. \aom{53}{1951}{332}.}
  \ref{Weil2}{Weil,A.\ma{168}{1967}{}.}
  \ref{CandN}{Chandrasekharan,K. and Narasimhan,R. \aom{74}{1961}{1}.}
  \ref{Rankin3}{Rankin,R.A. {} {} ().}
  \ref{Berndt6}{Berndt,B.C. {\it Trans.Edin.Math.Soc}.}
  \ref{Elizalde}{Elizalde,E. {\it Ten Physical Applications of Spectral
  Zeta Function Theory}, \break (Springer, Berlin, 1995).}
  \ref{Allen}{Allen,B., Folacci,A. and Gibbons,G.W. \pl{189}{1987}{304}.}
  \ref{Krazer}{Krazer}
  \ref{Elizalde3}{Elizalde,E. {\it J.Comp.and Appl. Math.} {\bf 118}
  (2000) 125.}
  \ref{Elizalde2}{Elizalde,E., Odintsov.S.D, Romeo, A. and Bytsenko,
  A.A and
  Zerbini,S.
  {\it Zeta function regularisation}, (World Scientific, Singapore,
  1994).}
  \ref{Eisenstein}{Eisenstein}
  \ref{Hecke}{Hecke,E. \ma{112}{1936}{664}.}
  \ref{Hecke2}{Hecke,E. \ma{112}{1918}{398}.}
  \ref{Terras}{Terras,A. {\it Harmonic analysis on Symmetric Spaces} (Springer,
  New York, 1985).}
  \ref{BandG}{Bateman,P.T. and Grosswald,E. {\it Acta Arith.} {\bf 9}
  (1964) 365.}
  \ref{Deuring}{Deuring,M. \aom{38}{1937}{585}.}
  \ref{Guinand}{Guinand.}
  \ref{Guinand2}{Guinand.}
  \ref{Minak}{Minakshisundaram.}
  \ref{Mordell}{Mordell,J. \prs{}{}{}.}
  \ref{GandZ}{Glasser,M.L. and Zucker, {}.}
  \ref{Landau2}{Landau,E. \jram{}{1903}{64}.}
  \ref{Kirsten1}{Kirsten,K. \jmp{35}{1994}{459}.}
  \ref{Sommer}{Sommer,J. {\it Vorlesungen \"uber Zahlentheorie}
  (1907,Teubner,Leipzig).
  French edition 1913 .}
  \ref{Reid}{Reid,L.W. {\it Theory of Algebraic Numbers},
  (1910,MacMillan,New York).}
  \ref{Milnor}{Milnor, J. {\it Is the Universe simply--connected?},
  IAS, Princeton, 1978.}
  \ref{Milnor2}{Milnor, J. \ajm{79}{1957}{623}.}
  \ref{Opechowski}{Opechowski,W. {\it Physica} {\bf 7} (1940) 552.}
  \ref{Bethe}{Bethe, H.A. \zfp{3}{1929}{133}.}
  \ref{LandL}{Landau, L.D. and Lishitz, E.M. {\it Quantum
  Mechanics} (Pergamon Press, London, 1958).}
  \ref{GPR}{Gibbons, G.W., Pope, C. and R\"omer, H., \np{157}{1979}{377}.}
  \ref{Jadhav}{Jadhav,S.P. PhD Thesis, University of Manchester 1990.}
  \ref{DandJ}{Dowker,J.S. and Jadhav, S. \prD{39}{1989}{1196}.}
  \ref{CandM}{Coxeter, H.S.M. and Moser, W.O.J. {\it Generators and
  relations of finite groups} (Springer. Berlin. 1957).}
  \ref{Coxeter2}{Coxeter, H.S.M. {\it Regular Complex Polytopes},
   (Cambridge University Press, \break Cambridge, 1975).}
  \ref{Coxeter}{Coxeter, H.S.M. {\it Regular Polytopes}.}
  \ref{Stiefel}{Stiefel, E., J.Research NBS {\bf 48} (1952) 424.}
  \ref{BandS}{Brink, D.M. and Satchler, G.R. {\it Angular momentum theory}.
  (Clarendon Press, Oxford. 1962.).}
  %\ref{Racah1}
  \ref{Rose}{Rose}
  \ref{Schwinger}{Schwinger, J. {\it On Angular Momentum} in {\it Quantum Theory of
  Angular Momentum} edited by Biedenharn,L.C. and van Dam, H.
  (Academic Press, N.Y. 1965).}
  \ref{Bromwich}{Bromwich, T.J.I'A. {\it Infinite Series},
  (Macmillan, 1947).}
  \ref{Ray}{Ray,D.B. \aim{4}{1970}{109}.}
  \ref{Ikeda}{Ikeda,A. {\it Kodai Math.J.} {\bf 18} (1995) 57.}
  \ref{Kennedy}{Kennedy,G. \prD{23}{1981}{2884}.}
  \ref{Ellis}{Ellis,G.F.R. {\it General Relativity} {\bf2} (1971) 7.}
  \ref{Dow8}{Dowker,J.S. \cqg{20}{2003}{L105}.}
  \ref{IandY}{Ikeda, A and Yamamoto, Y. \ojm {16}{1979}{447}.}
  \ref{BandI}{Bander,M. and Itzykson,C. \rmp{18}{1966}{2}.}
  \ref{Schulman}{Schulman, L.S. \pr{176}{1968}{1558}.}
  \ref{Bar1}{B\"ar,C. {\it Arch.d.Math.}{\bf 59} (1992) 65.}
  \ref{Bar2}{B\"ar,C. {\it Geom. and Func. Anal.} {\bf 6} (1996) 899.}
  \ref{Vilenkin}{Vilenkin, N.J. {\it Special functions},
  (Am.Math.Soc., Providence, 1968).}
  \ref{Talman}{Talman, J.D. {\it Special functions} (Benjamin,N.Y.,1968).}
  \ref{Miller}{Miller, W. {\it Symmetry groups and their applications}
  (Wiley, N.Y., 1972).}
  \ref{Dow3}{Dowker,J.S. \cmp{162}{1994}{633}.}
  \ref{Cheeger}{Cheeger, J. \jdg {18}{1983}{575}.}
  \ref{Cheeger2}{Cheeger, J. \aom {109}{1979}{259}.}
  \ref{Dow6}{Dowker,J.S. \jmp{30}{1989}{770}.}
  \ref{Dow20}{Dowker,J.S. \jmp{35}{1994}{6076}.}
  \ref{Dowjmp}{Dowker,J.S. \jmp{35}{1994}{4989}.}
  \ref{Dow21}{Dowker,J.S. {\it Heat kernels and polytopes} in {\it
   Heat Kernel Techniques and Quantum Gravity}, ed. by S.A.Fulling,
   Discourses in Mathematics and its Applications, No.4, Dept.
   Maths., Texas A\&M University, College Station, Texas, 1995.}
  \ref{Dow9}{Dowker,J.S. \jmp{42}{2001}{1501}.}
  \ref{Dow7}{Dowker,J.S. \jpa{25}{1992}{2641}.}
  \ref{Warner}{Warner.N.P. \prs{383}{1982}{379}.}
  \ref{Wolf}{Wolf, J.A. {\it Spaces of constant curvature},
  (McGraw--Hill,N.Y., 1967).}
  \ref{Meyer}{Meyer,B. \cjm{6}{1954}{135}.}
  \ref{BandB}{B\'erard,P. and Besson,G. {\it Ann. Inst. Four.} {\bf 30}
  (1980) 237.}
  \ref{PandM}{Polya,G. and Meyer,B. \cras{228}{1948}{28}.}
  \ref{Springer}{Springer, T.A. Lecture Notes in Math. vol 585 (Springer,
  Berlin,1977).}
  \ref{SeandT}{Threlfall, H. and Seifert, W. \ma{104}{1930}{1}.}
  \ref{Hopf}{Hopf,H. \ma{95}{1925}{313}. }
  \ref{Dow}{Dowker,J.S. \jpa{5}{1972}{936}.}
  \ref{LLL}{Lehoucq,R., Lachi\'eze-Rey,M. and Luminet, J.--P. {\it
  Astron.Astrophys.} {\bf 313} (1996) 339.}
  \ref{LaandL}{Lachi\'eze-Rey,M. and Luminet, J.--P.
  \prp{254}{1995}{135}.}
  \ref{Schwarzschild}{Schwarzschild, K., {\it Vierteljahrschrift der
  Ast.Ges.} {\bf 35} (1900) 337.}
  \ref{Starkman}{Starkman,G.D. \cqg{15}{1998}{2529}.}
  \ref{LWUGL}{Lehoucq,R., Weeks,J.R., Uzan,J.P., Gausman, E. and
  Luminet, J.--P. \cqg{19}{2002}{4683}.}
  \ref{Dow10}{Dowker,J.S. \prD{28}{1983}{3013}.}
  \ref{BandD}{Banach, R. and Dowker, J.S. \jpa{12}{1979}{2527}.}
  \ref{Jadhav2}{Jadhav,S. \prD{43}{1991}{2656}.}
  \ref{Gilkey}{Gilkey,P.B. {\it Invariance theory,the heat equation and
  the Atiyah--Singer Index theorem} (CRC Press, Boca Raton, 1994).}
  \ref{BandY}{Berndt,B.C. and Yeap,B.P. {\it Adv. Appl. Math.}
  {\bf29} (2002) 358.}
  \ref{HandR}{Hanson,A.J. and R\"omer,H. \pl{80B}{1978}{58}.}
  \ref{Hill}{Hill,M.J.M. {\it Trans.Camb.Phil.Soc.} {\bf 13} (1883) 36.}
  \ref{Cayley}{Cayley,A. {\it Quart.Math.J.} {\bf 7} (1866) 304.}
  \ref{Seade}{Seade,J.A. {\it Anal.Inst.Mat.Univ.Nac.Aut\'on
  M\'exico} {\bf 21} (1981) 129.}
  \ref{CM}{Cisneros--Molina,J.L. {\it Geom.Dedicata} {\bf84} (2001)
  \ref{Goette1}{Goette,S. \jram {526} {2000} 181.}
  207.}
  \ref{NandO}{Nash,C. and O'Connor,D--J, \jmp {36}{1995}{1462}.}
  \ref{Dows}{Dowker,J.S. \aop{71}{1972}{577}; Dowker,J.S. and Pettengill,D.F.
  \jpa{7}{1974}{1527}; J.S.Dowker in {\it Quantum Gravity}, edited by
  S. C. Christensen (Hilger,Bristol,1984)}
  \ref{Jadhav2}{Jadhav,S.P. \prD{43}{1991}{2656}.}
  \ref{Dow11}{Dowker,J.S. \cqg{21}{2004}4247.}
  \ref{Dow12}{Dowker,J.S. \cqg{21}{2004}4977.}
  \ref{Dow13}{Dowker,J.S. \jpa{38}{2005}1049.}
  \ref{Zagier}{Zagier,D. \ma{202}{1973}{149}}
  \ref{RandG}{Rademacher, H. and Grosswald,E. {\it Dedekind Sums},
  (Carus, MAA, 1972).}
  \ref{Berndt7}{Berndt,B, \aim{23}{1977}{285}.}
  \ref{HKMM}{Harvey,J.A., Kutasov,D., Martinec,E.J. and Moore,G.
  {\it Localised Tachyons and RG Flows}, hep-th/0111154.}
  \ref{Beck}{Beck,M., {\it Dedekind Cotangent Sums}, {\it Acta Arithmetica}
  {\bf 109} (2003) 109-139 ; math.NT/0112077.}
  \ref{McInnes}{McInnes,B. {\it APS instability and the topology of the brane
  world}, hep-th/0401035.}
  \ref{BHS}{Brevik,I, Herikstad,R. and Skriudalen,S. {\it Entropy Bound for the
  TM Electromagnetic Field in the Half Einstein Universe}; hep-th/0508123.}
  \ref{BandO}{Brevik,I. and Owe,C.  \prD{55}{4689}{1997}.}
  \ref{Kenn}{Kennedy,G. Thesis. University of Manchester 1978.}
  \ref{KandU}{Kennedy,G. and Unwin S. \jpa{12}{L253}{1980}.}
  \ref{BandO1}{Bayin,S.S.and Ozcan,M.
  \prD{48}{2806}{1993}; \prD{49}{5313}{1994}.}
  \ref{Chang}{Chang, P., {\it Quantum Field Theory on Regular Polytopes}.
   Thesis. University of Manchester, 1993.}
  \ref{Barnesa}{Barnes,E.W. {\it Trans. Camb. Phil. Soc.} {\bf 19} (1903) 374.}
  \ref{Barnesb}{Barnes,E.W. {\it Trans. Camb. Phil. Soc.}
  {\bf 19} (1903) 426.}
  \ref{Stanley1}{Stanley,R.P. \joa {49Hilf}{1977}{134}.}
  \ref{Stanley}{Stanley,R.P. \bams {1}{1979}{475}.}
  \ref{Hurley}{Hurley,A.C. \pcps {47}{1951}{51}.}
  \ref{IandK}{Iwasaki,I. and Katase,K. {\it Proc.Japan Acad. Ser} {\bf A55}
  (1979) 141.}
  \ref{IandT}{Ikeda,A. and Taniguchi,Y. {\it Osaka J. Math.} {\bf 15} (1978)
  515.}
  \ref{GandM}{Gallot,S. and Meyer,D. \jmpa{54}{1975}{259}.}
  \ref{Flatto}{Flatto,L. {\it Enseign. Math.} {\bf 24} (1978) 237.}
  \ref{OandT}{Orlik,P and Terao,H. {\it Arrangements of Hyperplanes},
  Grundlehren der Math. Wiss. {\bf 300}, (Springer--Verlag, 1992).}
  \ref{Shepler}{Shepler,A.V. \joa{220}{1999}{314}.}
  \ref{SandT}{Solomon,L. and Terao,H. \cmh {73}{1998}{237}.}
  \ref{Vass}{Vassilevich, D.V. \plb {348}{1995}39.}
  \ref{Vass2}{Vassilevich, D.V. \jmp {36}{1995}3174.}
  \ref{CandH}{Camporesi,R. and Higuchi,A. {\it J.Geom. and Physics}
  {\bf 15} (1994) 57.}
  \ref{Solomon2}{Solomon,L. \tams{113}{1964}{274}.}
  \ref{Solomon}{Solomon,L. {\it Nagoya Math. J.} {\bf 22} (1963) 57.}
  \ref{Obukhov}{Obukhov,Yu.N. \pl{109B}{1982}{195}.}
  \ref{BGH}{Bernasconi,F., Graf,G.M. and Hasler,D. {\it The heat kernel
  expansion for the electromagnetic field in a cavity}; math-ph/0302035.}
  \ref{Baltes}{Baltes,H.P. \prA {6}{1972}{2252}.}
  \ref{BaandH}{Baltes.H.P and Hilf,E.R. {\it Spectra of Finite Systems}
  (Bibliographisches Institut, Mannheim, 1976).}
  \ref{Ray}{Ray,D.B. \aim{4}{1970}{109}.}
  \ref{Hirzebruch}{Hirzebruch,F. {\it Topological methods in algebraic
  geometry} (Springer-- Verlag,\break  Berlin, 1978). }
  \ref{BBG}{Bla\v{z}i\'c,N., Bokan,N. and Gilkey, P.B. {\it Ind.J.Pure and
  Appl.Math.} {\bf 23} (1992) 103.}
  \ref{WandWi}{Weck,N. and Witsch,K.J. {\it Math.Meth.Appl.Sci.} {\bf 17}
  (1994) 1017.}
  \ref{Norlund}{N\"orlund,N.E. \am{43}{1922}{121}.}
  \ref{Duff}{Duff,G.F.D. \aom{56}{1952}{115}.}
  \ref{DandS}{Duff,G.F.D. and Spencer,D.C. \aom{45}{1951}{128}.}
  \ref{BGM}{Berger, M., Gauduchon, P. and Mazet, E. {\it Lect.Notes.Math.}
  {\bf 194} (1971) 1. }
  \ref{Patodi}{Patodi,V.K. \jdg{5}{1971}{233}.}
  \ref{GandS}{G\"unther,P. and Schimming,R. \jdg{12}{1977}{599}.}
  \ref{MandS}{McKean,H.P. and Singer,I.M. \jdg{1}{1967}{43}.}
  \ref{Conner}{Conner,P.E. {\it Mem.Am.Math.Soc.} {\bf 20} (1956).}
  \ref{Gilkey2}{Gilkey,P.B. \aim {15}{1975}{334}.}
  \ref{MandP}{Moss,I.G. and Poletti,S.J. \plb{333}{1994}{326}.}
  \ref{BKD}{Bordag,M., Kirsten,K. and Dowker,J.S. \cmp{182}{1996}{371}.}
  \ref{RandO}{Rubin,M.A. and Ordonez,C. \jmp{25}{1984}{2888}.}
  \ref{BaandD}{Balian,R. and Duplantier,B. \aop {112}{1978}{165}.}
  \ref{Kennedy2}{Kennedy,G. \aop{138}{1982}{353}.}
  \ref{DandKi2}{Dowker,J.S. and Kirsten, K. {\it Analysis and Appl.}
 {\bf 3} (2005) 45.}
  \ref{Dow40}{Dowker,J.S. \cqg{23}{2006}{1}.}
  \ref{BandHe}{Br\"uning,J. and Heintze,E. {\it Duke Math.J.} {\bf 51} (1984)
   959.}
  \ref{Dowl}{Dowker,J.S. {\it Functional determinants on M\"obius corners};
    Proceedings, `Quantum field theory under
    the influence of external conditions', 111-121,Leipzig 1995.}
  \ref{Dowqg}{Dowker,J.S. in {\it Quantum Gravity}, edited by
  S. C. Christensen (Hilger, Bristol, 1984).}
  \ref{Dowit}{Dowker,J.S. \jpa{11}{1978}{347}.}
  \ref{Kane}{Kane,R. {\it Reflection Groups and Invariant Theory} (Springer,
  New York, 2001).}
  \ref{Sturmfels}{Sturmfels,B. {\it Algorithms in Invariant Theory}
  (Springer, Vienna, 1993).}
  \ref{Bourbaki}{Bourbaki,N. {\it Groupes et Alg\`ebres de Lie}  Chap.III, IV
  (Hermann, Paris, 1968).}
  \ref{SandTy}{Schwarz,A.S. and Tyupkin, Yu.S. \np{242}{1984}{436}.}
  \ref{Reuter}{Reuter,M. \prD{37}{1988}{1456}.}
  \ref{EGH}{Eguchi,T. Gilkey,P.B. and Hanson,A.J. \prp{66}{1980}{213}.}
  \ref{DandCh}{Dowker,J.S. and Chang,Peter, \prD{46}{1992}{3458}.}
  \ref{APS}{Atiyah M., Patodi and Singer,I.\mpcps{77}{1975}{43}.}
  \ref{Donnelly}{Donnelly.H. {\it Indiana U. Math.J.} {\bf 27} (1978) 889.}
  \ref{Katase}{Katase,K. {\it Proc.Jap.Acad.} {\bf 57} (1981) 233.}
  \ref{Gilkey3}{Gilkey,P.B.\invm{76}{1984}{309}.}
  \ref{Degeratu}{Degeratu.A. {\it Eta--Invariants and Molien Series for
  Unimodular Groups}, Thesis MIT, 2001.}
  \ref{Seeley}{Seeley,R. \ijmp {A\bf18}{2003}{2197}.}
  \ref{Seeley2}{Seeley,R. .}
  \ref{melrose}{Melrose}
  \ref{berard}{B\'erard,P.}
  \ref{gromes}{Gromes,D.}
  \ref{Ivrii}{Ivrii}
  \ref{DandW}{Douglas,R.G. and Wojciekowski,K.P. \cmp{142}{1991}{139}.}
  \ref{Dai}{Dai,X. \tams{354}{2001}{107}.}
  \ref{Kuznecov}{Kuznecov}
  \ref{DandG}{Duistermaat and Guillemin.}
  \ref{PTL}{Pham The Lai}
\end{putreferences}

\bye